\documentclass[10pt]{elsart}

\usepackage{amsmath}
\usepackage{amsopn}
\usepackage{amssymb}
\usepackage[all]{xy}
\usepackage{graphicx}
\newtheorem{mythm}{Theorem}[subsection]
\newtheorem{mycor}[mythm]{Corollary}
\newtheorem{mylem}[mythm]{Lemma}
\newtheorem{myprop}[mythm]{Proposition}
\newtheorem{myrem}[mythm]{Remark}
\newtheorem{mydefn}[mythm]{Definition}

\newcommand{\abs}[1]{\lvert #1 \rvert}

\newcommand{\br}[1]{\overline{#1}}

\newcommand{\td}[1]{\widetilde{#1}}

\newcommand{\bra}[1]{\langle #1 \rangle}

\newcommand{\ZZ}{\mathbb{Z}}

\newcommand{\QQ}{\mathbb{Q}}
\newcommand{\WW}{\mathbb{W}}
\newcommand{\FF}{\mathbb{F}}
\newcommand{\GG}{\mathbb{G}}
\newcommand{\MS}{\mathbb{S}}

\newcommand{\HH}{\mathbb{H}}
\newcommand{\calM}{\mathcal{M}}
\newcommand{\calO}{\mathcal{O}}

\newcommand{\calC}{\mathcal{C}}

\numberwithin{equation}{subsection}
\numberwithin{figure}{section}
\DeclareMathOperator{\ext}{Ext}
\DeclareMathOperator{\Hom}{Hom}

\DeclareMathOperator{\spec}{spec}
\DeclareMathOperator{\map}{Map}
\DeclareMathOperator{\aut}{Aut}
\DeclareMathOperator{\stack}{stack}
\DeclareMathOperator{\Frob}{Frob}
\DeclareMathOperator{\End}{End}

\begin{document}

\begin{frontmatter}

\title{A modular description of the $K(2)$-local sphere at the prime $3$}
\author{Mark Behrens\thanksref{NSF}}
\address{Department of Mathematics, MIT, Cambridge, MA
  02139, USA}
\thanks[NSF]{The author is partially supported by the NSF.}

\begin{abstract}
Using degree $N$ isogenies of elliptic curves, we produce a spectrum $Q(N)$.  
This spectrum is built out of spectra related to $tmf$.
At $p=3$ we show that the $K(2)$-local sphere is built out of $Q(2)$ and
its $K(2)$-local Spanier-Whitehead dual.  This gives a conceptual
reinterpretation a resolution of Goerss, Henn, Mahowald, and Rezk.

\end{abstract}


\end{frontmatter}

\tableofcontents
\setcounter{tocdepth}{1}

\section*{Introduction}

Goerss, Henn, Mahowald, and Rezk \cite{GHMR} have produced at the prime $3$ 
a tower of spectra of the form
\begin{multline*}
TMF \rightarrow TMF \vee \Sigma^8 E(2) \rightarrow \Sigma^8E(2) \vee
\Sigma^{40}E(2)
\\
\rightarrow \Sigma^{40}E(2) \vee \Sigma^{48}TMF \rightarrow \Sigma^{48}TMF.
\end{multline*}
Here, as is the case everywhere in this paper unless indicated
otherwise, we are working in the $K(2)$-local category, and everything has
been implicitly $K(2)$-localized.  
The authors of \cite{GHMR} show that the tower \emph{refines} 
to a resolution of the 
$K(2)$-local sphere at the prime
$3$.  This means that there exists a diagram
$$
\xymatrix{
L_{K(2)}S \ar[d] 
& X_1 \ar[l] \ar[d] 
& X_2 \ar[l] \ar[d]
& X_3 \ar[l] \ar[d] 
& X_4 \ar[l] \ar@{=}[d]
\\
TMF 
& \genfrac{}{}{0pt}{}{\Sigma^{-1}TMF \vee}{\Sigma^7E(2)}
& \genfrac{}{}{0pt}{}{\Sigma^6E(2) \vee}{\Sigma^{38}E(2)}
& \genfrac{}{}{0pt}{}{\Sigma^{37}E(2) \vee}{\Sigma^{45}TMF}
& \Sigma^{44}TMF
}
$$
such that each spectrum $X_i$ is the fiber of the $i$th vertical map,
and such that the connecting morphisms of these fiber sequences give the
maps in the tower.
This resolution helps organize the computation of Shimomura and Wang
of $\pi_*(L_{K(2)}S^0)$.  We shall refer to this as the GHMR resolution.

The spectrum $tmf$ is the connective 
spectrum of topological modular forms.  We write $TMF$ for the
non-connective spectrum
$tmf[\Delta^{-1}]$, but this distinction is irrelevant, since we have
$$ L_{K(2)}tmf \simeq L_{K(2)}TMF \simeq EO_2 = E_2^{hG_{24}}. $$
Here $E_2$ is Morava $E$-theory and $G_{24}$ is a maximal finite subgroup
of the extended Morava stabilizer group $\GG_2$ at the prime $3$.  The
spectrum $EO_2$ exists by the 
Hopkins-Miller theorem \cite{RezkHopkinsMiller}, and
is $E_\infty$ by the machinery of Goerss and Hopkins \cite{GoerssHopkins}.
The equivalence $L_{K(2)}TMF \simeq EO_2$ is discussed in greater detail in
Remark~\ref{rmk:K(2)localTMF}.

\begin{pf*}{Remark.}
It has become standard to let $E_2$ denote the Landweber exact cohomology
theory corresponding to the Lubin-Tate deformation of the 
Honda height $2$ formal group $F_2$, defined over $\FF_{p^2}$, with
$p$-series
$$ [p]_{F_2} = x^{p^2}. $$
However, for the $3$-primary applications in this paper, we shall let $E_2$
denote the spectrum corresponding to the Lubin-Tate deformation of the
height $2$ formal group $F'_2$ over $\FF_{9}$ whose $3$-series is given by
$$ [3]_{F'_2} = -x^9. $$
Our reason for doing this is that $F'_2$ is the formal group of the unique
supersingular elliptic curve over $\FF_9$
(see Section~\ref{sec:aut}, in particular Remark~\ref{rmk:Galoisaction}).
\end{pf*}

The GHMR resolution is a generalization to chromatic
level $2$ of the $J$-spectrum fiber sequence
$$ L_{K(1)}S^0 \rightarrow KO_p \xrightarrow{\psi^N-1} KO_p. $$
That this gives the $K(1)$-local sphere was first observed in unpublished
work of Adams and Baird, and an independent verification of was
provided by Ravenel \cite{Ravenel84}, \cite{Bousfield}.
Here $N$ is a topological generator of $\ZZ_p^\times$.

There is a spectrum $TMF_0(2)$ which is an analog of $TMF$ for the congruence
subgroup $\Gamma_0(2) < SL_2(\ZZ)$.  We have
equivalences
$$ L_{K(2)}TMF_0(2) \simeq E_2^{hD_8} \simeq L_{K(2)}(E(2) \vee \Sigma^8
E(2))$$
where $D_8 < \GG_2$ is a dihedral group of order $8$.  After
inverting $6$, the
spectrum $TMF_0(2)$ coincides with the elliptic cohomology
theory $Ell$ of Landweber, Ravenel, and Stong
\cite{LandweberRavenelStong}. 

The GHMR tower may be made less efficient, to
take the form
\begin{multline*}
TMF \rightarrow TMF \vee TMF_0(2) \rightarrow TMF_0(2) \vee
\Sigma^{48}TMF_0(2)
\\
\rightarrow \Sigma^{48}TMF_0(2) \vee \Sigma^{48}TMF
\rightarrow \Sigma^{48}TMF  
\end{multline*}
by letting the extra copies of $E(2)$ kill themselves pairwise.  It is this 
resolution that we shall endeavor to reinterpret.

The GHMR resolution was produced using the work of Devinatz and Hopkins
\cite{DevinatzHopkins},
where the $K(2)$-local sphere is identified as a homotopy fixed point spectrum
$$ E_2^{h\GG_2} \simeq L_{K(2)}S^0. $$ 
The authors of \cite{GHMR} produced a resolution of the trivial
$\GG_2$-module by permutation modules.  
The authors then realize their resolution to
produce the tower.

Mahowald and Rezk wanted a modular description of the GHMR resolution. 
The motivation is that the Adams operation $\psi^N$
corresponds to the $N^{th}$ power isogeny on the multiplicative group $\GG_m$.
By replacing the multiplicative group with an elliptic curve, one can
instead consider certain degree $N$ isogenies of elliptic curves.  
Mahowald and Rezk studied the
corresponding map
$$ TMF \xrightarrow{"\psi^N-1"} TMF_0(N) $$
at the primes $2$ and $3$ \cite{MahowaldRezk}, \cite{Mahowald}.
This setup is complicated by the fact that there is a whole moduli
space of elliptic curves, and elliptic curves do not support unique degree 
$N$ isogenies.  The map that Mahowald and Rezk studied is a 
projection of the first map of
the GHMR resolution.  Ando, Hopkins, and
Strickland have also extensively studied operations on elliptic cohomology
theories arising from isogenies of elliptic curves \cite{AndoEn},
\cite{AndoElliptic}, \cite{AndoHopkinsStrickland}.

The aim of this paper is expand on the work of Mahowald and Rezk to produce
a tower
$$ TMF \rightarrow TMF_0(N) \vee TMF \rightarrow TMF_0(N) $$
which refines to a resolution of an $E_\infty$ ring spectrum $Q(N)$.
Here $N$ is prime to $p$.
This tower arises conceptually from certain degree $N$ 
isogenies of elliptic curves.
Rezk had considered a similar setup from the point of view of flags of
subgroups of an elliptic curve.
For a spectrum $X$, let $D_{K(2)}(X)$ denote the Spanier-Whitehead dual in
the $K(2)$-local category
$$ D_{K(2)}(X) = F(X,L_{K(2)}S). $$
We conjecture that for appropriate $N$, the spectrum 
$Q(N)$ is half of the $K(2)$-local sphere. 

\begin{conj}\label{conj:cofiber}
Let $p$ be greater than $2$.
If $N$ is chosen to be a 
topological generator
of $\ZZ_p^\times$, then the natural sequence
$$ D_{K(2)}Q(N) \xrightarrow{D\eta} L_{K(2)} S \xrightarrow{\eta} Q(N) $$
is a cofiber sequence.  Here $\eta$ is the unit of the $K(2)$-local ring
spectrum $Q(N)$.

If $p = 2$, and $N$ is a topological generator of an index $2$ subgroup of
$\ZZ_2^\times$, then there is an index $2$ subgroup $\td{\GG}_2$ of $\GG_2$
and a cofiber sequence
$$ D_{K(2)}Q(N) \xrightarrow{D\eta} E_2^{h\td{\GG}_2} 
\xrightarrow{\eta} Q(N). $$
\end{conj}

Note that the group $\ZZ_2^\times$ is not topologically cyclic, but it 
does have
topologically cyclic index $2$ subgroups.

The main result of this paper is to verify Conjecture~\ref{conj:cofiber} 
at the prime $3$, with
$N=2$ (Theorem~\ref{thm:cofiber}).  The author also knows
Conjecture~\ref{conj:cofiber} to hold in the case $p=5$ and $N=2$.  The
author has no evidence for the invariance of Conjecture~\ref{conj:cofiber}
under different choices of $N$.

This approach to understanding the $K(2)$-local sphere has several advantages.  The maps building $Q(N)$ are quite
computable using well known formulas for Weierstrass curves and have 
very natural descriptions in terms of isogenies of elliptic
curves.  
The cofiber sequence of Conjecture~\ref{conj:cofiber} explains the
self-duality of the GHMR resolution, and an identification of $D_{K(2)}TMF$
explains the appearance of the suspension $\Sigma^{48}$.  The very difficult
computations of Shimomura and Wang of the homotopy of $L_{K(2)} S$ should be
verified independently using our decomposition.  As a side-effect of our
work we also reproduce the short tower (Proposition~\ref{prop:Sbarlagrangian})
\begin{equation}\label{eq:Sbartower}
TMF \rightarrow \Sigma^8 E(2) \rightarrow \Sigma^{40} E(2) \rightarrow
\Sigma^{48} TMF
\end{equation}
which refines to the spectrum $\br{S} = E^{h\GG^1_2}$ 
given in \cite{GHMR}.  Here
$\GG^1_2$ is the kernel of the reduced norm
$$ \GG_2 \rightarrow \ZZ_3. $$

We are not able to describe the connecting map
$$ Q(N) \rightarrow \Sigma D_{K(2)}Q(N) $$
of Conjecture~\ref{conj:cofiber}, nor are we able to describe the middle
map of the short tower (\ref{eq:Sbartower}).

Our approach to proving Theorem~\ref{thm:cofiber} is
computational, and this is the reason for our specialization to the prime
$3$ and $N=2$.
It would be nice to have a conceptual and elegant proof of
Conjecture~\ref{conj:cofiber} for all primes and all $N$.  
The $2$-primary applications with $N=3$ should be rewarding.  There is
essentially nothing to be gained computationally from
Conjecture~\ref{conj:cofiber} if $p \ge 5$.

This paper is organized into two parts.  A detailed 
outline of the content is given
at the beginning of each part.
In Part~1 we give a construction of $Q(N)$, and quickly
specialize to the case $N=2$ and $p=3$.  In this case we give a computation
of the $V(1)$-homology 
groups $V(1)_*Q(2)$ where $V(1)$ is the Smith-Toda complex.
In Part~2 we prove Conjecture~\ref{conj:cofiber} in the
case $N=2$ and $p=3$ (Theorem~\ref{thm:cofiber}).  
The proof uses the $V(1)$ homology
computations given in Part~1 as well as the computation of
the $V(1)$ homology of $L_{K(2)}S^0$ --- and thus
does not generalize to arbitrary $N$ and $p$.

\emph{Everything in this paper is implicitly $K(2)$-localized unless
specifically specified otherwise.  Sections for which this convention does
not hold are indicated as such in their beginnings.}

We highlight a few aspects of this paper which may be of independent
interest.  The duality decomposition of Conjecture~\ref{conj:cofiber} may
be interpreted as a Lagrangian decomposition for a hyperbolic duality
pairing on the $K(2)$-local sphere $L_{K(2)}S$.  
The abstract framework of duality
pairings and Lagrangians in a triangulated symmetric monoidal category is
given in Section~\ref{sec:pairings}.

In Section~\ref{sec:Sbar}, we show that for arbitrary $n$, at an arbitrary
prime $p$, the canonical
pairing on $L_{K(n)}S$ is hyperbolic, and $\br{S}=E_n^{h\GG^1_n}$ 
is a Lagrangian.  We identify the $K(n)$-local Spanier-Whitehead dual of
$\br{S}$ by proving there is an equivalence
$$ D_{K(n)}(\br{S}) \simeq \Sigma^{-1} \br{S}. $$

Finally, for computational reasons, we need a cellular decomposition of
$\br{S}$ in the $K(n)$-local category.  This task is relegated to
Appendix~A, where we prove that there is a ($K(n)$-local) 
cellular decomposition
$$ \br{S} \simeq S^0 \cup_{\zeta} e^0 \cup_{\zeta} e^0 \cup_{\zeta} \cdots. $$

The author would like to thank Charles Rezk and Mark Mahowald, 
who shared their ideas and preliminary notes so
generously.  Thanks also go to Mike Hopkins for suggesting the 
duality mechanism which is the main result of this paper, 
Paul Goerss for explaining
$E_n$-homology operations, Sharon Hollander and Tilman Bauer, for explaining
various aspects of
stacks, and Nasko Karamanov for sharing his computational knowledge.  
Daniel Davis provided valuable assistance with homotopy fixed point
spectra, in particular with Theorem~\ref{thm:FEHEH} and
Lemma~\ref{lem:telnull}.
Many ideas in this paper are culled from various works of
Matthew Ando and Neil Strickland.  Discussions with 
Peter May, Haynes Miller, and Doug
Ravenel during
the course of this project were very helpful.  
The initial version of this paper dealt with the Galois action on the
Morava stabilizer group improperly.  The author was alerted to this problem
by Hans-Werner Henn, and Mike Hopkins' help in sorting it out was
invaluable.
Thanks also go to the
referee for pointing out that the hypotheses of Proposition~\ref{prop:MRGH}
needed modification.

\section*{Part 1: The spectrum $Q(N)$}
\addtocounter{section}{1}
\addcontentsline{toc}{part}{Part 1: The spectrum $Q(N)$ \hfill}

This part is organized as follows.
In Section~\ref{sec:simpstack} we produce a simplicial stack based on
certain degree $N$ isogenies.  

In Section~\ref{sec:realize} we topologically 
realize our simplicial stack to a cosimplicial spectrum whose totalization
is our spectrum $Q(N)$.  We give two approaches to this realization
problem.  The first is based on a sheaf of $E_\infty$ ring spectra produced
by Hopkins and his collaborators.  
The second, specialized to $p=3$ and $N=2$, uses the
Goerss-Hopkins-Miller theory and a supersingular elliptic curve.

In Section~\ref{sec:level2} we give a description of the ring of
$\Gamma_0(2)$
modular forms in terms of Weierstrass equations.  The results of this
section are probably well known.

Our computations of the homotopy of $Q(2)$ are based on the ANSS.  
In Section~\ref{sec:ANcomplex}, we describe a chain complex which computes
the Adams-Novikov $E_2$ term for $Q(2)$.  This ANSS $E_2$-term is the
hypercohomology of the simplicial stack we
constructed in Section~\ref{sec:simpstack}.

Section~\ref{sec:mapscomp} is devoted to explicitly computing 
the maps on
rings of modular forms induced by the face maps of our simplicial stack.
From these maps we get the differentials in the chain complex of
Section~\ref{sec:ANcomplex}.  These explicit formulas are also
necessary to complete the construction of $Q(2)$ given in
Section~\ref{sec:sheafrealize}.

In Section~\ref{sec:aut} we describe explicitly the automorphism groups of
a supersingular elliptic curve $C$ and its formal group $C^\wedge$.
In Section~\ref{sec:Morava} we use this description to 
describe the effects that our maps have on
$E_2$-homology.  
These formulas will be used in many technical lemmas in
Part~2, and also may be used in the
construction of $Q(2)$ given in Section~\ref{sec:superrealize}.

Our approach to Theorem~\ref{thm:cofiber} is computational and based on
$V(1)$-homology computations.  In Section~\ref{sec:Q(2)V(1)}, we compute
the $V(1)$-homology of $Q(2)$.  This computation is compared to
computations of Goerss, Henn, Mahowald, and Shimomura.

\subsection{A simplicial stack}\label{sec:simpstack}

Let $\mathcal{M}$ be the moduli stack of non-singular elliptic curves over
$\ZZ_{(p)}$.
Let $\omega$ be the line bundle of invariant differentials over
$\mathcal{M}$, so that the sections of $\omega^{\otimes k}$ give weight $k$
weakly modular forms.
$$ H^0(\mathcal{M}, \omega^{\otimes k}) = MF_k $$
$MF_*$ is the ring of weakly modular forms over $\ZZ_{(p)}$ 
\cite{Deligne}.
$$ MF_* = \ZZ_{(p)}[c_4, c_6, \Delta^{\pm1}]/(c_4^3-c_6^2 = 1728\Delta ) $$
We must invert $\Delta$ because we are taking sections over the complement
of the singular locus of the moduli space of generalized elliptic curves.
We shall refer to the ring of modular forms (without inverting the
discriminant) as $mf_*$.

Given an elliptic curve $C$, we shall mean by a $\Gamma_0(N)$ structure 
a chosen discrete subgroup of $C$ which is isomorphic to $C_N$, 
the cyclic group
of order $N$, but we do \emph{not} fix the isomorphism.  $\Gamma_0(N)$ 
structures
are in one to one correspondence
with degree $N$ isogenies with cyclic kernels.  Given a $\Gamma_0(N)$
structure $H$ of $C$, one has an isogeny
$$ \phi_H: C \rightarrow C/H $$
and given such an isogeny, one recovers the $\Gamma_0(N)$ structure by taking
the kernel.  
We define $\calM_0(N)$ to be the moduli stack of non-singular elliptic curves
over $\ZZ_{(p)}$ with $\Gamma_0(N)$ structure.  We shall always be considering
$N$ prime to $p$.

Moduli stacks of elliptic curves with additional structure have been widely
studied by arithmetic algebraic geometers.  See, for example,
\cite{DeligneRapoport}, \cite{KatzMazur}.

Given an integer $N$, every elliptic curve has an $N^{th}$ power
endomorphism 
$$ [N] : C \rightarrow C $$
given (for positive $N$) by
$$ [N](P) = \underbrace{P + \cdots + P}_{N}. $$
It is an isogeny of degree $N^2$.

Given any degree $N$ isogeny $\phi$, there exists a 
dual isogeny $\widehat{\phi}$.  The dual isogeny has the property that 
$$ \widehat{\phi} \circ \phi = [N]. $$
See, for instance, \cite[III.6]{Silverman}.

We shall now describe a simplicial stack $\mathcal{M_\bullet}$.  
Actually $\calM_\bullet$ is a semi-simplicial stack
since we do not use degeneracies.  It 
has only $0$, $1$, and $2$ simplices.  
We will describe this simplicial object first in conceptual terms, and then in 
precise stack theoretic language.
The simplices will be given as follows:
\begin{itemize}
\item The $0$-simplices are elliptic curves.
\item The $1$-simplices are certain isogenies.  There will be two types:
  \begin{itemize}
  \item Degree $N$ isogenies $\phi_H$ with cyclic kernel $H$.  These are in
  one-to-one correspondence with elliptic curves $C$ with a $\Gamma_0(N)$ 
  structure $H$.
  \item The endomorphisms $[N]$.  These are in one-to-one correspondence with
  elliptic curves, since every elliptic curve possesses a unique such
  endomorphism.
  \end{itemize}
\item The $2$-simplices 
correspond to relations of the form $\widehat{\phi}_H \circ \phi_H =
[N]$.  These relations are in one-to-one correspondence with 
elliptic curves with
a $\Gamma_0(N)$ structure $H$, 
since there is one such relation for every isogeny
$\phi_H$.
\end{itemize}

The semi-simplicial stack
$\mathcal{M_\bullet}$ is a diagram
$$
\xymatrix@C+1em{
\calM & 
\calM \coprod \calM_0(N) \ar@<-1ex>[l]|-{d_0} \ar@<1ex>[l]|-{d_1} & 
\calM_0(N) \ar@<-1.5ex>[l]|-{d_0} \ar[l]|-{d_1} \ar@<1.5ex>[l]|-{d_2}
}
$$
in the category of stacks.
Given a $\ZZ_{(p)}$-algebra $R$, the $R$-points
of $\calM$ is a groupoid of elliptic curves $C$ over $R$, and
the $R$-points of $\calM_0(N)$ is a groupoid of pairs 
$(C,H)$ of elliptic curves $C$ over
$R$ with $\Gamma_0(N)$
structure $H \le C(R)$. The simplicial stack $\calM_\bullet$ should be
viewed on the level of $R$-points as consisting of moduli of diagrams of
the following form.
$$
\begin{array}{ccccc}
\left\{ C \right\} 
& 
\begin{array}{c}
\xleftarrow{d_0} \\
\xleftarrow{d_1}
\end{array}
&
\begin{array}{c}
\left\{ (C,H) \xrightarrow{\phi_H} C/H  \right\} \\
\coprod \\
\left\{ C \xrightarrow{[N]} C \right\}
\end{array}
&
\begin{array}{c}
\xleftarrow{d_0} \\
\xleftarrow{d_1} \\
\xleftarrow{d_2}
\end{array}
&
\left\{
\xymatrix@C-1em{
(C,H) \ar[rr]^{[N]} \ar[rd]_{\phi_H} && C  \\
& C/H \ar[ru]_{\widehat{\phi}_H}
}
\right\}
\end{array}
$$
The face maps are what one might expect from taking the nerve of a
category.  
To give them we must define some maps of stacks.  These maps of stacks 
are given on
the level of $R$-points by the formulas below.
\begin{itemize}
\item Forget $\Gamma_0(N)$ structure:
\begin{gather*}
\phi_f: \calM_0(N) \rightarrow \calM \\
(C,H) \mapsto C 
\end{gather*}
\item Quotient out by level $\Gamma_0(N)$ structure:
\begin{gather*}
\phi_q : \calM_0(N) \rightarrow \calM \\
(C,H) \mapsto C/H 
\end{gather*}
\item Codomain of the $N^{th}$ power endomorphism:
\begin{gather*}
\psi_{[N]} : \calM \rightarrow \calM \\
C \mapsto C/C[N]
\end{gather*}
\item Dual $\Gamma_0(N)$ structure 
(the dual $\Gamma_0(N)$ structure $\widehat{H}$ 
is defined such that $\phi_{\widehat{H}} =
\widehat{\phi}_H$):
\begin{gather*}
\psi_d : \calM(N) \rightarrow \calM(N) \\
(C,H) \mapsto (C/H, \widehat{H})
\end{gather*}
\end{itemize}
We comment that there is a canonical isomorphism
$$ r_C : C/C[N] \cong C. $$
This isomorphism is \emph{not} an equality, and this distinction is important
in the context of stacks.  
We may use this isomorphism to define $\psi_{[N]}$ on $\calM_0(N)$ by
\begin{gather*}
\psi_{[N]}: \calM_0(N) \rightarrow \calM_0(N) \\
(C,H) \mapsto (C/C[N],r_C^{-1}(H)).
\end{gather*}
We note that we have the following 
relations
\begin{gather}
\phi_q(C) = \phi_f(\psi_d(C)) \label{eq:rel1} \\
\psi_d \circ \psi_d = \psi_{[N]}. \label{eq:rel2}
\end{gather}
The first relation indicates that if we have defined $\phi_f$ and $\psi_d$,
then $\phi_q$ is automatically determined.

The careful reader will be bothered that the maps $\phi_f$, $\psi_{[N]}$,
and $\psi_d$ are not defined above in a precise manner.  In
Section~\ref{sec:mapscomp}, we give explicit formulas for $N=2$ and $p=3$ 
for these maps on the level of
Hopf algebroids.  The stackifications of these Hopf algebroids 
are $\calM$ and $\calM_0(2)$, so we get induced maps of stacks.  We can then
check relations~\ref{eq:rel1} and \ref{eq:rel2} explicitly.  Since we are
only concerned with the case $p=3$, $N=2$ in this paper, we choose to not
elaborate further on the general case.

We use these maps of stacks to give the face maps of $\calM_\bullet$.
The face maps $d_i: \calM \coprod \calM_0(N) \rightarrow \calM$ are defined
by
\begin{gather*}
d_0 = \psi_{[N]} \coprod \phi_q \\
d_1 = Id_{\calM} \coprod \phi_f
\end{gather*}
The face maps $d_i: \calM_0(N) \rightarrow \calM \coprod \calM_0(N)$ are
defined by
\begin{gather*}
d_0 = \psi_d \\
d_1 = \phi_f \\
d_2 = Id_{\calM_0(N)}
\end{gather*}
The simplicial identities are verified by relations~\ref{eq:rel1} and
\ref{eq:rel2}.

\subsection{Realizing the simplicial stack}\label{sec:realize}

We wish to topologically realize the simplicial stack $\calM_\bullet$ as a
cosimplicial $E_\infty$ ring spectrum $Q(N)^\bullet$.  Actually,
$Q(N)^\bullet$ is a semi-cosimplicial $E_\infty$ ring spectrum, since we have
no codegeneracies.
$$
\xymatrix@C+1em{
TMF \ar@<-1ex>[r]|-{d_0} \ar@<1ex>[r]|-{d_1} & 
TMF \times TMF_0(N) 
\ar@<-1.5ex>[r]|-{d_0} \ar[r]|-{d_1} \ar@<1.5ex>[r]|-{d_2} & 
TMF_0(N) 
}
$$
We then define the spectrum $Q(N)$ by
$$ Q(N) = Tot(Q(N)^\bullet). $$
The spectrum $Q(N)$ is an $E_\infty$ ring spectrum since the coface maps 
are maps of $E_\infty$ ring spectra.

We give two
approaches to constructing $Q(N)^\bullet$.  
The first is based on a construction of $tmf$ due to
Hopkins and his collaborators 
as the global sections of a sheaf of $E_\infty$ ring spectra.  Since
this work is unpublished, we give an alternative construction for $p=3$ and
$N=2$ using the
Goerss-Hopkins refinement of the Hopkins-Miller theorem.  The latter
approach is sufficient for our purposes since we are working
$K(2)$-locally.

\subsubsection{Sheaf theoretic construction}\label{sec:sheafrealize}

Hopkins and his collaborators have constructed a 
sheaf $\calO_{ell}$ of $E_\infty$ ring spectra on 
$\calM$ in the \'etale topology.  
Paul Goerss has written a useful survey of this point of view \cite{Goerss}.
The sheaf $\calO_{ell}$ has the property that if $C$ is
any elliptic curve over $R$ which is \'etale over $\calM$, then the
homotopy groups of the sections over $C$ are given by
$$ \pi_{2k}(\calO_{ell}(C)) \cong \omega^{\otimes k}_C $$
where $\omega_C$ is the free $R$-module of holomorphic 1-forms on $C$.
One recovers the spectra $TMF$ and $TMF_0(N)$ as the sections of this sheaf
over the appropriate moduli stacks.
\begin{align*}
TMF & = \calO_{ell}(\calM) \\
TMF_0(N) & = \calO_{ell}(\calM_0(N))
\end{align*}
Since $\calO_{ell}$ is a sheaf in the \'etale topology,
we need the map
$$ \phi_f: \calM_0(N) \rightarrow \calM $$
to be \'etale to take sections over $\calM_0(N)$.  
This is why we insist that $p$ does not divide $N$.

Our requirement that $p$ does not divide $N$ also implies that the maps 
\begin{gather*}
\psi_{[N]}: \calM_0(N) \rightarrow \calM_0(N) \\
\psi_{[N]}: \calM \rightarrow \calM \\
\end{gather*}
are isomorphisms.  Relation~\ref{eq:rel2} implies that 
$$ \psi_d : \calM_0(N) \rightarrow \calM_0(N) $$
is an isomorphism.  We have already indicated that the map $\phi_f$ is
\'etale.  Relation~\ref{eq:rel1} indicates that $\phi_q$ is \'etale.  We
therefore get induced maps
\begin{align*}
\phi_f^* : & TMF = \calO_{ell}(\calM) \rightarrow \calO_{ell}(\calM_0(N)) =
TMF_0(N) \\
\phi_q^* : & TMF = \calO_{ell}(\calM) \rightarrow \calO_{ell}(\calM_0(N)) =
TMF_0(N) \\
\psi_{[N]}^* : & TMF = \calO_{ell}(\calM) \rightarrow
\calO_{ell}(\calM) = TMF \\
\psi_d^* : & TMF_0(N) = \calO_{ell}(\calM_0(N)) \rightarrow
\calO_{ell}(\calM_0(N)) = TMF_0(N).
\end{align*}
These maps induce the coface maps $d_i$ of $Q(N)^\bullet$.

\subsubsection{Supersingular construction}\label{sec:superrealize}

We now give an alternative construction of $Q(2)^\bullet$ at $p=3$ that uses
computations of endomorphisms of the supersingular elliptic curve of
Section~\ref{sec:aut}.
This approach to $TMF$ mirrors the approach in \cite{HopkinsMahowald}.
The Goerss-Hopkins-Miller 
theorem \cite{RezkHopkinsMiller}, \cite{GoerssHopkins} 
states that there is a contravariant functor
$$ E : \mathcal{FGL}_n \rightarrow \text{\it $E_\infty$ ring spectra} $$
where $\mathcal{FGL}_n$ is the category of pairs $(k,F)$ where $k$ is a
perfect field of characteristic $p$ and $F$ is a formal group law 
of height $n$ over $k$.  Consider the supersingular elliptic curve
$$ C: y^2 = x^3 - x $$
over $\FF_9$.  The Weierstrass curve $C$ has a canonical coordinate,
and gives rise to a formal group law $C^\wedge$ of height $2$.
The curve $C$ has a $\Gamma_0(2)$ structure $H$ generated by the point
$(x,y)=(0,0)$ (see Section~\ref{sec:level2}).  
In Section~\ref{sec:Morava}, we produce an explicit degree
$2$ isogeny
$$ \psi_H: C \rightarrow C $$
with $\ker \psi_H = H$.  We have extended automorphism groups (including 
automorphisms of the ground field $\FF_9$) of elliptic curves 
(respectively, elliptic curves with $\Gamma_0(2)$
structures) given by
\begin{align*}
\aut_{/\FF_3}(C) & = G_{24} \\
\aut_{/\FF_3}(C,H) & = D_8.
\end{align*}
Here $G_{24}$ and $D_8$ correspond to certain subgroups of $\GG_2$.  These
subgroups are described explicitly in Sections~\ref{sec:aut} and 
\ref{sec:Morava}.
The map $\psi_H$ is invariant under $D_8$, since every element of
$D_8$ fixes the $\Gamma_0(2)$ structure $H$.

The value of the Goerss-Hopkins-Miller functor $E$ on
$(\FF_9,C^\wedge)$ is the version of Morava $E$-theory 
that we will be using.
$$ E_{(\FF_9,C^\wedge)} = E_2 $$
The map $\psi_H$ induces an isomorphism of formal group laws
$$ \widehat{\psi}_H: C^\wedge \rightarrow C^\wedge $$
invariant under the action of $D_8$, and 
the Goerss-Hopkins-Miller theorem implies that it induces a map of
$E_\infty$ ring spectra
$$ \psi_d^* : TMF_0(2) = E_2^{hD_8} \rightarrow E_2^{hD_8} = TMF_0(2). $$
We verify directly in
Section~\ref{sec:Morava} that $(\psi_d^*)^2 = \psi_{[2]}^*$, where
$\psi_{[2]}^*$
is the Adams operation
corresponding to the action of the element $2 \in \WW_{\FF_9}^\times
\subset \MS_2$ of the Morava stabilizer group.  Define
$\phi_f^*$ to be the restriction
$$ \phi_f^* = Res^{G_{24}}_{D_8}: TMF = E_2^{hG_{24}} \rightarrow
E_2^{hD_8}. $$
We get $\phi_q^*$ for free by defining it to be the composite
$$ \phi_q^* : TMF \xrightarrow{\phi_f^*} TMF_0(2) \xrightarrow{\psi_d^*}
TMF_0(2). $$
We now can define the coface maps $d_i$ in terms of $\phi_f^*$,
$\psi_{[2]}^*$, $\psi_d^*$, and $\phi_q^*$ as before, and we get our
cosimplicial spectrum $Q(2)^\bullet$.

The idea is that $K(2)$-locally, we only need to consider sections of
$\calO_{ell}$ in a formal neighborhood of the supersingular locus, and that
is precisely what the Goerss-Hopkins-Miller theorem is implicitly doing at
chromatic level $2$.

\subsection{$\Gamma_0(2)$ modular forms}\label{sec:level2}

\subsubsection{Moduli description of $mf_0(2)$}
In order to compute the maps in the resolution for $p=3$ and $N=2$ we need
to know something about $\Gamma_0(2)$ modular forms.  
Most of the computations in this section have been carried out by Mark Mahowald
and Charles Rezk, or are just in the literature on
modular forms.
It is well known (see for instance, the appendix of
\cite{HirzebruchBergerJung}) that the
ring of $\Gamma_0(2)$ modular forms is given (with $2$ inverted) by
$$ mf_0(2)_*[1/2] = \ZZ[1/2][\delta, \epsilon] $$
where the weights are given by
\begin{align*}
\abs{\delta} & = 2 \\
\abs{\epsilon} & = 4.
\end{align*}
The discriminant is given by
$$ \Delta = -64\epsilon^2(\epsilon+\delta^2). $$
Denote the ring with $\Delta$ inverted by
$$ MF_0(2)_*[1/2] = mf_0(2)_*[1/2,\Delta^{-1}]. $$
Let $\calM_0(2)$ be the moduli stack of 
non-singular elliptic curves with $\Gamma_0(2)$
structure.
The ring $MF_0(2)$ 
should be interpreted as sections of tensor powers of the line
bundle $\omega$ over this stack.
$$
MF_0(2)_k  = H^0(\calM_0(2), \omega^{\otimes k})
$$

We can eliminate our use of the line bundle $\omega$ if we add the data of
a non-zero tangent vector to the structures on elliptic curves 
we are taking moduli of.
Namely, let $\calM^1_0(2)$ be the moduli stack of 
non-singular elliptic curves with the
data of a
$\Gamma_0(2)$ structure and a non-zero tangent vector at the identity.  
Then we have
$$
MF_0(2)_* = H^0(\calM^1_0(2), \calO)
$$
where $\calO$ is the structure sheaf of $\calM^1_0(2)$.
There is an action of the multiplicative group $\GG_m$ on $\calM_0^1(2)$; the
action is given by the multiplication action on the chosen tangent vector.
Under this action, the sections above break up into a direct sum as the
weights of the $\GG_m$ representations.  These coincide with the weights of
the modular forms.

\subsubsection{Calculation of $MF_0(2)_*[1/2]$ using Weierstrass
curves}\label{sec:mf(2)[1/2]}

We shall give an alternate description of the generators of $MF_0(2)_*[1/2]$
using Weierstrass equations.  For the rest of this section we will be
implicitly working with $2$ inverted.  Under this condition, every
Weierstrass curve is isomorphic to one of the form \cite{Silverman}
$$ C_\mathbf{b} : y^2 = 4x^3 + b_2x^2 + 2b_4x + b_6. $$
This curve is nonsingular if the discriminant
$$ \Delta = -27b_6^2 + (9b_2b_4 - \frac{1}{4}b_2^3)b_6 - 8b_4^3 + 
\frac{1}{4}b_2^2b_4^2 $$
is a unit.  We shall only consider non-singular curves.
We implicitly think of these curves as coming with the data of the tangent
vector $\partial/\partial z$ where $z = x/y$.  The only transformations
which leave this form of equation and tangent vector invariant are those of
the form
$$ \chi_r: x \mapsto x+r, y \mapsto y. $$
A $\Gamma_0(2)$ structure is the choice of a point of exact order $2$ in
$C_\mathbf{b}$.  These points coincide with the points $(x,y)$ 
of the curve $C_{\mathbf{b}}$
with $y = 0$.  Thus a $\Gamma_0(2)$ structure consists of a chosen root
$e_1$ of the
cubic $4x^3 + b_2x^2 + 2b_4x + b_6$.  The curves $C_\mathbf{b}$ with this
additional data may be written in the form
$$ C_{\mathbf{\gamma}, e_1} : y^2 = 4(x-e_1)(x^2 + \gamma_2x + \gamma_4) $$
where we have
\begin{align*}
b_2 & = 4(\gamma_2-e_1) \\
b_4 & = 2(\gamma_4-e_1\gamma_2) \\
b_6 & = -4\gamma_4e_1.
\end{align*}
Given the root $e_1$ and the coefficients $b_2$, $b_4$, $b_6$, the
coefficients $\gamma_2$ and $\gamma_4$ may be expressed as
\begin{align*}
\gamma_2 &= \frac{b_2+e_1}{4} \\
\gamma_4 &= \frac{4b_4+e_1b_2+e_1^2}{8}.
\end{align*}
The data of the $\Gamma_0(2)$ structure allows us to remove the automorphism
$\chi_r$.  We simply shift the root $e_1$ to be $0$ by the transformation
$\chi_{e_1}$.  This puts $C_{\mathbf{\gamma},e_1}$ into the canonical form
$$ C_{\mathbf{q}}: y^2 = 4x(x^2+q_2x+q_4). $$
The curve $C_\mathbf{q}$ is regarded as implicitly carrying the data of
the
tangent vector $\partial/\partial z$ and the 
$\Gamma_0(2)$ structure generated by
the point $(x,y) = (0,0)$.  After the transformation $\chi_{e_1}$, one has
\begin{align}
q_2 & = 2e_1+\gamma_2 \label{eq:qgamma1} \\
q_4 & = e_1^2 + \gamma_2e_1 + \gamma_4.  \label{eq:qgamma2}
\end{align}
The discriminant of the curve $C_\mathbf{q}$ is given by
$$ \Delta = q_4^2(16q_2^2-64q_4). $$
There are no non-trivial automorphisms of the curves $C_\mathbf{q}$ 
preserving the $\Gamma_0(2)$ structure and the tangent vector, 
thus we have shown
that the stack $\calM^1_0(2)$ is an affine scheme.
$$ \calM^1_0(2) = \spec(\ZZ[1/2, q_2, q_4, \Delta^{-1}]) $$
The ring $MF_0(2)_*[1/2]$ is just the ring of functions.
$$ MF_0(2)_*[1/2] = \ZZ[1/2, q_2, q_4, \Delta^{-1}] $$

\subsubsection{Relation to $\delta$ and $\epsilon$}

We now relate our generators $q_2$ and $q_4$ to the classical generators
$\delta$ and $\epsilon$ that appear in the Jacobi quartic.  Assume we are
working over an algebraically closed field $k$ of characteristic not equal
to $2$ or $3$.  An elliptic curve can be expressed by a
Weierstrass equation
$$ C_\mathbf{e} : y^2 = 4(x-e_1)(x-e_2)(x-e_3) $$
with $e_1+e_2+e_3 = 0$.  (This last condition is equivalent to the
Weierstrass equation taking the form
$ y^2 = 4x^3 - g_2x-g_3$.)
We can take $e_1$ to be the specified $\Gamma_0(2)$ structure.  In
\cite{HirzebruchBergerJung}, the quantities $\delta$ and $\epsilon$ are
then given by
\begin{align*}
\delta & = - \frac{3}{2} e_1 \\
\epsilon & = (e_1-e_2)(e_1-e_3).
\end{align*}
By multiplying out the factors containing $e_2$ and $e_3$, we see that
$C_\mathbf{e}$, with $\Gamma_0(2)$ structure $e_1$, is given by the curve
$C_{\mathbf{\gamma},e_1}$ with
\begin{align*}
\gamma_2 & = -e_2-e_3 \\
\gamma_4 & = e_2e_3.
\end{align*}
Using Equations~\ref{eq:qgamma1} and \ref{eq:qgamma2}, and the relation
$e_1+e_2+e_3 = 0$, we see that
\begin{align*}
q_2 & = -2\delta \\
q_4 & = \epsilon.
\end{align*}
Thus our generators $q_2$ and $q_4$ of $mf_0(2)[1/2]$ are, up to scaler
multiple, identical to the classical generators $\delta$ and $\epsilon$.

\subsection{The Adams-Novikov $E_2$-term of
$Q(2)$}\label{sec:ANcomplex}

There are several approaches to computing 
the homotopy groups of the spectrum $Q(2)$.  One could
compute the maps that $\phi_q$, $\phi_f$, $\psi_d$, and $\psi_{[2]}$
induce on the homotopy groups of $TMF$ and $TMF_0(2)$, and then use the
Bousfield-Kan spectral sequence for $Tot(Q(2)^\bullet)$, but it is actually
easier to compute the $E_2$-term of the Adams-Novikov spectral sequence
(ANSS) for $Q(2)$ and then add in the
Adams-Novikov differentials using the differentials in the ANSS for $TMF$.
The discussion of the ANSS for $TMF$ in this section 
follows \cite{Rezk} and \cite{HopkinsMahowald}.  Everything in this section
is implicitly $3$-local, though most of 
what we say holds if we simply invert $2$.

\subsubsection{The ANSS $E_2$ term for $TMF$.}

We explain how to compute the ANSS $E_2$-term for $TMF$ using an \'etale
cover of $\calM^1$.  This material has appeared elsewhere, see
\cite{HopkinsMahowald}, \cite{Rezk}, and \cite{Bauer}. 
As mentioned in Section~\ref{sec:mf(2)[1/2]}, over a
$\ZZ_{(3)}$-algebra 
every elliptic curve is isomorphic to one of the
form 
$$ C_{\mathbf{b}} : y^2 = 4x^3+b_2x^2 +2b_4 x + b_6. $$
The automorphisms of this elliptic curve which preserve the tangent vector
are of the form
$$ \chi_r: x \mapsto x+r. $$
Consider the elliptic curve Hopf algebroid
\begin{align*}
B & = \ZZ_{(3)}[b_2,b_4, b_6, \Delta^{-1}] \\
\Gamma_B & = B[r]
\end{align*}
which represents the groupoid of such elliptic curves and isomorphisms.
The right unit is determined by how the coefficients of the Weierstrass
equation for $C_{\mathbf{b}}$ transform.  
Since every isomorphism class of elliptic 
curve is represented, the stackification (in the flat topology) of 
the pre-stack associated to the
Hopf algebroid $(B,\Gamma_B)$ is $\calM^1$.

The Hopf algebroid $(B, \Gamma_B)$ 
suffers from the drawback that the natural map
$$ \spec (B) \rightarrow \calM^1 $$
is not \'etale.  
We instead consider the curves
$$ C_{\mathbf{b}} : y^2 = 4x(x^2+q_2x+q_4). $$
The automorphisms of these curves which preserve the tangent vector are
those of the form $x \mapsto x+r$ where $r$ has the property that  
$$ r^3+ q_2 r^2 + q_4r = 0. $$
This groupoid of elliptic curves is represented by the Hopf algebroid
$(\br{B},\Gamma_{\br{B}})$ given by
\begin{align*}
\br{B} & = \ZZ_{(3)}[q_2,q_4, \Delta^{-1}] \\
\Gamma_{\br{B}} & = \br{B}[r]/(r^3+ q_2 r^2 + q_4r).
\end{align*}

\begin{myprop}\label{prop:etalecover}
The map
$$ f_{\mathbf{q}} : \spec(\br{B}) \rightarrow \calM^1$$ 
which classifies the curve
$C_{\mathbf{q}}$
is an \'etale cover.
\end{myprop}

\begin{pf*}{Proof.}
We first point out that the natural map of stackifications
$$ \stack(\br{B}, \Gamma_{\br{B}}) \rightarrow \stack(B, \Gamma_B) $$
is an equivalence of stacks.  This follows since given a curve $C$ of the form
$C_{\mathbf{b}}$ over a $\ZZ_{(p)}$-algebra $R$, there is a finite faithfully
flat extension 
$R' = R[t]/(4t^3+b_2t^2 + 2b_4t + b_6)$
of $R$ so that over $R'$, $C$ is isomorphic to a curve of the form
$C_{\mathbf{q}}$ under the transformation $x \mapsto x+t$.

Therefore we simply must check that the natural map
$$ \spec(\br{B}) \rightarrow \stack(\br{B}, \Gamma_{\br{B}}) $$
is an \'etale cover.  It suffices to check that the pullback
$$
\xymatrix{
\spec(\Gamma_{\br{B}}) \ar[r] \ar[d] &
\spec(\br{B}) \ar[d] 
\\
\spec(\br{B}) \ar[r] &
\stack(\br{B}, \Gamma_{\br{B}})
}
$$
is an \'etale cover, that is, that $\Gamma_{\br{B}}$ is an \'etale
extension.  This follows from the fact that since $\Delta$ is invertible,
$q_4$ is invertible,
and therefore the derivative of $f(r) = r^3 + q_2r^2+q_4r$ is
nonzero modulo all maximal ideals in $\br{B}$.
\qed \end{pf*}

In particular, $f_{\mathbf{q}}$ is a flat cover.  Thus we have the following
corollary.

\begin{mycor}
The ANSS $E_2$-term for $TMF$ may be computed as
$$ H^*(\calM^1, \calO) \cong \ext_{\Gamma_{\br{B}}}(\br{B}, \br{B}) \cong
H^*(C^*(\Gamma_{\br{B}})) $$
where $C^*(\Gamma_{\br{B}})$ is the cobar complex of the Hopf algebroid
$(\br{B}, \Gamma_{\br{B}})$.
\end{mycor}

Since the stack $\calM^1_0(2)$ is an affine scheme after
inverting $2$,
the ANSS for $TMF_0(2)$ is concentrated in the zero line.  The ANSS collapses
to give
$$ \pi_{2k}(TMF_0(2)) = MF(2)_k. $$

\subsubsection{The ANSS $E_2$-term for $Q(2)$}

The ANSS for $Q(2)$ may be obtained by totalizing a cosimplicial
Adams-Novikov resolution for $Q(2)$.  Therefore, its $E_2$-term is the
hypercohomology of the simplicial stack $\calM_\bullet$.  For
convenience we choose to work with the simplicial stack $\calM^1_\bullet$ 
where all of the
instances of $\calM$ and $\calM_0(2)$ are replaced by $\calM^1$ and
$\calM^1_0(2)$, respectively.  The resulting spectral sequence takes the form
$$ \HH^{s,t}(\calM^1_\bullet, \calO) \Rightarrow \pi_{2t-s}(Q(2)).  $$
The $E_2$ term may be computed via a hypercohomology spectral sequence
\begin{equation}\label{eq:ANE2SS}
H^*_{cosimp}(H^*(\calM^1_\bullet, \calO)) \Rightarrow
\HH^*(\calM^1_\bullet, \calO).
\end{equation}
Here $H^*_{cosimp}$ is the cohomology of the resulting cosimplicial abelian
group.

The hypercohomology group $\HH^*(\calM^1_\bullet,\calO)$ 
is the cohomology of the totalization of a
double cochain complex $C^{*,*}(Q(2))$ whose horizontal differentials are
given by 
$$
\xymatrix@C+.49em@R-1em{
C^*(\Gamma_{\br{B}}) \ar[r]^-{d_0-d_1} \ar@{=}[d] &
MF_0(2) \oplus \br{C}^*(\Gamma_{\br{B}}) \ar[r]^-{d_0-d_1+d_2} \ar@{=}[d] &
MF_0(2) \ar[r] \ar@{=}[d] &
0 \ar[r] \ar@{=}[d] &
\,
\\
C^{0,*}(Q(2)) &
C^{1,*}(Q(2)) &
C^{2,*}(Q(2)) &
C^{3,*}(Q(2)) 
}$$
and whose vertical differentials are given by the differentials of the cobar
complex.
Here $MF_0(2)$ should be regarded as a
cochain complex concentrated in degree zero.  The complex
$\br{C}^*(\Gamma_{\br{B}})$ is the cobar complex where the differentials have 
been given the opposite sign.  The coface maps $d_i$ are lifts of the maps
induced by the simplicial face maps of $\calM^1_\bullet$.  In
Section~\ref{sec:mapscomp} we shall compute these maps explicitly.  

We
shall let the cochain complex $C^*(Q(2))$ be the totalization of the double complex 
$C^{*,*}(Q(2))$.  Then we
have
$$ \HH^*(\calM^1_\bullet,\calO) = H^*(C^*(Q(2))) $$
and the hypercohomology spectral sequence is simply the spectral sequence of
the double complex.

\subsection{Computation of the maps}\label{sec:mapscomp}

In this section we compute the effects of the
maps $\phi_q$, $\phi_f$, $\psi_d$, and
$\psi_{[2]}$ on the appropriate rings of modular forms.  
Actually, since
$\calM^1$ is not a scheme, we will lift the maps involving
$\calM^1$ to the prestack associated to the Hopf algebroid 
$(\br{B}, \Gamma_{\br{B}})$.  Our computations
of these maps will show, as a side-effect, that $\phi_f$ is \'etale, and
$\psi_{[2]}$ is an isomorphism.  This was what was required to
topologically realize the maps on the spectra of sections of $\calO_{ell}$
in Section~\ref{sec:sheafrealize}.  We shall also see that we have the
relation
$$ \psi_d^2 = \psi_{[2]}. $$
This appeared as relation~\ref{eq:rel2}, and it is required for the
simplicial identities to hold in the semi-simplicial stack $\calM_\bullet$.

\subsubsection{The map $\phi_f$}\label{sec:phi_f}

The stack $\calM^1_0(2)$ is the affine scheme $\spec(\br{B})$.  There are no
nontrivial automorphisms.  Forgetting the $\Gamma_0(2)$ structure generated by 
$(x,y) = (0,0)$ is the same as allowing the automorphisms of the curve
$C_{\mathbf{q}}$ which do not preserve this level $\Gamma_0(2)$ structure.  
Thus the map
$$ \phi_f : \calM^1_0(2) \rightarrow \calM^1 $$
is induced by the map of Hopf algebroids
$$ (\br{B}, \Gamma_{\br{B}}) \rightarrow (\br{B}, \br{B}) $$
which is the identity on $\br{B}$ and maps $r$ to zero.
Proposition~\ref{prop:etalecover} implies that the map $\phi_f$ is \'etale.

\subsubsection{The map $\psi_{[2]}$}\label{sec:psi_[2]}

The map $\psi_{[2]}: \calM^1 \rightarrow \calM^1$ 
takes an elliptic curve $C$ and returns the quotient $C/C[2]$ where $C[2]$ is
the subgroup of all points of order $2$ of $C$.
The quotient map $C \rightarrow C/C[2]$ is equivalent to the 
second power map $[2]$ on the elliptic curve.  The second power map on the
curve $C_{\mathbf{q}}$ could be
regarded as an endomorphism of the curve, but the tangent vector
$\partial/\partial z$ changes, where $z = x/y$.  Therefore we shall instead
regard the map $[2]$ as being a map from $C_\mathbf{q}$ to $C_\mathbf{q'}$
where the tuple $\mathbf{q'}$ is given by
\begin{align*}
q_2' & = 2^2 q_2 \\
q_4' & = 2^4 q_4.
\end{align*}
We also must adjust the automorphism $x \mapsto x+r$ to $x' \mapsto x'+r'$
where 
$$ r' = 2^2 r. $$
Therefore, the map on Hopf algebroids
$$ \psi_{[2]}^* : (\br{B}, \Gamma_{\br{B}}) \rightarrow (\br{B},
\Gamma_{\br{B}}) $$
is given by
\begin{align*}
\psi^*_{[2]}(q_2) & = 2^2 q_2 \\
\psi^*_{[2]}(q_4) & = 2^4 q_4 \\
\psi^*_{[2]}(r) & = 2^2 r.
\end{align*}
The map
$$ \psi_{[2]}^* : \calM^1_0(2) \rightarrow \calM^1_0(2) $$
is the map above restricted to $\br{B}$.  Each of these instances of
$\psi_{[2]}$ is an isomorphism since we can explicitly define the inverse
$\psi_{[1/2]}$ on the Hopf algebroids, since we are working $3$-locally.

\subsubsection{The map $\psi_d$}\label{sec:psi_d}

The dual isogeny map $\psi_d: \calM^1_0(2) \rightarrow \calM^1_0(2)$ takes a
curve with $\Gamma_0(2)$ structure $(C,H)$,
identifies $H$ as the kernel of a degree $2$ isogeny 
$$ \phi_H : C \rightarrow C/H, $$
and returns the elliptic curve with $\Gamma_0(2)$ structure 
$(C/H, \widehat{H})$
corresponding to the
dual isogeny
$$ \widehat{\phi}_{H}: C/H \rightarrow C. $$

In our case we are given the curve $C_\mathbf{q}$ with $\Gamma_0(2)$
structure $H$ generated by $e_1 = (x,y) = (0,0)$.  
We wish to find a Weierstrass equation for the quotient 
curve $C_{\mathbf{q}}/H$.
Given an elliptic curve $C$, a Weierstrass
equation for $C$ is generated by choosing a function $x_1$ with a pole of
order $2$ at the identity and a function $y_1$ with a pole of order $3$ at
the identity.  We need to find such functions $x_1$ and $y_1$ for
$C/H$.  

Define the functions $x_1$ and $y_1$ on $C_{\mathbf{q}}/H$ as follows.  
Given a point $P \in C_{\mathbf{q}}$, define 
\begin{align*}
x_1(P) & = x(P)+x_t(P) \\
y_1(P) & = y(P) + y_t(P).
\end{align*}
Here the functions $x_t$ and $y_t$ are defined by
\begin{align*}
x_t(P) &= x(P + e_1) \\
y_t(P) &= y(P + e_1)
\end{align*}
where $e_1$ is the generator of the $\Gamma_0(2)$ structure of $C_{\mathbf{q}}$ 
and $+$
denotes addition in the elliptic curve.  It is immediate that the
functions $x_1$ and $y_1$ descend to the quotient $C/H$.

One may use the formulas for
the group law of a Weierstrass curve given
in \cite{Silverman} to obtain
\begin{align*}
x_t & = \frac{q_4}{x} \\
y_t & = -\frac{q_4y}{x^2}.
\end{align*}
We warn the reader that the formulas of \cite{Silverman} must be modified
slightly since those formulas are for an elliptic curve in Weierstrass form 
$$ C_\mathbf{a}: y^2+a_1xy+a_3y = x^3+a_2x^2+a_4x+a_6 $$
whereas our curves
$$ C_\mathbf{q}: y^2 = 4x(x^2+q_2x+q_3) $$
have a coefficient of $4$ instead of $1$ in front of the $x^3$.

One then checks that the new coordinate functions $x_1$ and $y_1$ satisfy
the Weierstrass equation
\begin{equation}\label{eq:Weierquotient}
y_1^2 = 4x_1^3 + 4q_2x_1^2-16q_4x_1-16q_2q_4.
\end{equation}
The canonical tangent vector associated to the Weierstrass
equation~\ref{eq:Weierquotient} is the image of the tangent vector of
$C_{\mathbf{q}}$.
The right-hand side of equation~\ref{eq:Weierquotient} has a canonical root 
$x = -q_2$.
The curve $C_\mathbf{q}/H$ with this $\Gamma_0(2)$ structure $e_1 = -q_2$
corresponds to the dual isogeny to $C_\mathbf{q} \rightarrow
C_\mathbf{q}/H$.  We must put the Weierstrass
equation~\ref{eq:Weierquotient} into the form of $C_\mathbf{q'}$ for some
tuple $\mathbf{q'}$.  This is accomplished, as seen in
Section~\ref{sec:mf(2)[1/2]}, by application of the transformation
$x \mapsto x-q_2$.  The curve transforms to
$$
C_\mathbf{q'} : y^2 = 4x(x^2-2q_2x-4q_4+q_2^2).
$$
We conclude that the induced map
$$ \psi_d^*: \br{B} \rightarrow
\br{B} $$
is given by
\begin{align*}
\psi_d^*(q_2) & = -2q_2 \\
\psi_d^*(q_4) & = -4q_4+q_2^2.
\end{align*}
One sees immediately that relation~\ref{eq:rel2}
$$ (\psi_d)^2 = \psi_{[2]} : \calM^1_0(2) \rightarrow \calM^1_0(2) $$
is satisfied.  Since $\psi_{[2]}$ is an isomorphism, it follows that
$\psi_d$ is as well.

\subsubsection{The map $\phi_q$}\label{sec:phi_q}

Relation~\ref{eq:rel1} forces us do define $\phi_q$ to be the composite
$$ \phi_q : \calM^1_0(2) \xrightarrow{\psi_d} \calM^1_0(2) \xrightarrow{\phi_f}
\calM^1. $$
Therefore, on Hopf algebroids, the map 
$$ \phi_q^* : (\br{B}, \Gamma_{\br{B}}) \rightarrow (\br{B}, \br{B}) $$
is given by
\begin{align*}
\phi_q^*(q_2) & = -2q_2 \\
\phi_q^*(q_4) & = -4q_4 + q_2^2 \\
\phi_q^*(r) & = 0.
\end{align*}

\subsection{Extended automorphism groups}\label{sec:aut}

In this section we explain how to extend both the automorphism
group of an elliptic curve, and its associated formal group, by the Galois
group.  We then explicitly determine the extended automorphism group of the
unique supersingular elliptic curve over $\FF_9$, and the extended
automorphism group of its formal group.

In characteristic $p$, there are many different morphisms which all deserve
to be called the Frobenius morphism, and this can lead to confusion.  
We first recall some standard
material to clarify which one of these versions of the Frobenius we want 
to employ.

Let $q$ be a power of $p$.
Let $X = X' \times_{\spec(\FF_q)} \spec(\FF_{q^n})$ be a (formal) 
scheme over $\FF_{q^n}$ which is obtained from a (formal) scheme $X'$
over $\FF_q$ by base change.  
Let $\pi$ be the projection morphism
$$ \pi: X \rightarrow  \spec(\FF_{q^n}). $$
Let $\Frob_q: \FF_{q^n} \rightarrow \FF_{q^n}$ be the $q$th power Frobenius
automorphism, the generator of the Galois group $Gal(\FF_{q^n}/\FF_q)$.
Associated to $\Frob_q$ are several different  
Frobenius morphisms on
$X$:
\begin{itemize}
\item $\Frob_q$: the $q$th Frobenius on the base. \\
\item $\Frob^{rel}_q$: the $q$th relative Frobenius. \\
\item $\Frob^{tot}_q$: the $q$th total Frobenius.
\end{itemize}
These are given by the following diagram.
$$
\xymatrix@C+2em{
X \ar[ddr]_{\pi} \ar[rrd]^{\Frob^{tot}_q}
\ar[dr]|{\Frob^{rel}_q}  
\\
& X \ar[r]_{\Frob_q} \ar[d]_\pi
& X \ar[d]^\pi
\\
& \spec{\FF_{q^n}} \ar[r]_{\Frob_q} & \spec{\FF_{q^n}}
}
$$
The square in the above diagram is a pullback square.  If $X' = \spec(A)$
is affine, then $X$ is the affine scheme $\spec (A \otimes_{\FF_q}
\FF_{q^n})$, and $\Frob_q$ (respectively $\Frob_q^{rel}$) is the $q$th 
Frobenius
on $\FF_{q^n}$ (respectively on $A$), while $\Frob^{tot}_q$ is the Frobenius
on $A \otimes_{\FF_q} \FF_{q^n}$.

Let $\aut(X)$ denote the automorphism group of $X$, regarded as a (formal) 
scheme
over $\spec(\FF_{q^n})$  There is an induced action of $\Frob_q$ on
$\alpha \in \aut(X)$ given by the following pullback.
$$
\xymatrix@C+2em{
X \ar[r]^{\Frob^*_q \alpha} \ar[d]_{\Frob_q} &
X \ar[d]^{\Frob_q} 
\\
X \ar[r]_\alpha &
X
}
$$
This action of $\Frob_q$ gives rise to an action of 
$Gal(\FF_{q^n}/\FF_q)$ on $\aut(X)$. We define the
extended automorphism group $\aut_{/\FF_{q}}(X)$ to be the semidirect
product
$$ \aut_{/\FF_{q}}(X) = \aut(X) \rtimes Gal(\FF_{q^n}/\FF_q) $$
associated to this action.  The group $\aut_{/\FF_q}(X)$ consists of 
automorphisms of the (formal) scheme $X$ which do not cover the identity on
$\spec(\FF_{q^n})$.

At the prime $3$
there is one isomorphism class of supersingular elliptic curve over
$\FF_9$, admitting the Weierstrass presentation
$$ C : y^2 = x^3 - x. $$
For the remainder of this section we will focus our attention on the
computation of the endomorphism ring $\End(C)$, and the 
automorphism groups $\aut_{/\FF_3}(C)$ and 
$\aut_{/\FF_3}(C^\wedge)$.  We shall see that all of the endomorphisms
of $C$ defined over $\br{\FF}_3$ are defined over $\FF_9$.

Let $F = \Frob_3^{rel}$ be the relative Frobenius endomorphism of $C$.
The map $F$ is given by
$$ F: (x,y) \mapsto (x^3,y^3). $$
Define automorphisms $t, s \in \aut(C)$ by
\begin{align*}
t: & (x,y) \mapsto (\omega^4 x, \omega^6 y) \\
s: & (x,y) \mapsto (x + 1, y). 
\end{align*}
Here, $\omega \in \FF_9$ is a primitive $8$th root of unity.  Explicit
computation, using the formulas for the group law on $C$ \cite{Silverman}, 
demonstrates that we
have relations
\begin{align}
F^2 & = -3 \label{eq:phirel} \\  
t^2 & = -1 \label{eq:trel} \\
F t & = -tF \label{eq:phitrel} \\
s & = \frac{1}{2}(1 + F). \label{eq:sphirel}
\end{align}
From the definitions of the maps, we deduce the 
following additional relations.
\begin{align}
s^3 & = 1 \label{eq:srel} \\
s t& = t s^2 \label{eq:strel}
\end{align}

\begin{myprop}\label{prop:autC}
The automorphism group of $C$ is the group $G_{12} = C_3 \rtimes C_4$ 
of order $12$ generated 
by $s$ and $t$.
\end{myprop}

\begin{pf*}{Proof.}
Relations~(\ref{eq:srel}), (\ref{eq:trel}), and (\ref{eq:strel}) 
demonstrate that
$s$ and $t$ generate a subgroup $G_{12}$ of $\aut(C)$ of order $12$ 
isomorphic to
$C_3 \rtimes C_4$.  The group of automorphisms of $C$ defined over
$\br{\FF}_3$ is of order $12$ \cite[A.1.2]{Silverman}.  We therefore deduce
that all of the automorphisms of $C$ are defined over $\FF_9$, and that
they are contained in $G_{12}$.
\qed \end{pf*}

\begin{myprop}
The endomorphism ring of $C$ (over $\FF_9$) is a maximal order of 
the quaternion algebra
$$ \QQ \langle F, t \rangle / 
(F^2 = -3, t^2 = -1, F t = -tF). $$
A $\ZZ$-basis of the order $\mathcal{O} = \End(C)$ is given by
$$ \{ s, t, ts, tF\}. $$
\end{myprop}

\begin{pf*}{Proof.}
Since $C$ supersingular, the ring of endomorphisms defined 
over $\br{\FF}_3$ 
$\End_{\br{\FF}_3}(C)$ is a maximal order of a rational quaternion algebra
$D$ ramified at the places $p$ and $\infty$.  By Proposition~5.1 of 
\cite{Pizer}, $D$ 
admits the presentation
$$ D = \QQ \langle i,j \rangle / (i^2 = -3, j^2 = -1, ij = -ji). $$  By
Proposition~5.2 of \cite{Pizer}, there is a maximal order $\mathcal{O}$ of
$D$ with $\ZZ$-basis
$$ \{ \frac{1}{2} (1 + i), j, \frac{1}{2} (j + j i), j i \}.
$$
Let $R$ be the subring of $\End(C)$ generated by $t$, $s$, and $F$. 
Since $D$
is a division algebra, the algebra map
$$ D \rightarrow R \otimes \QQ $$
sending $i$ to $F$ and $j$ to $t$ is an isomorphism.  
We deduce that
under this isomorphism, the ring $R$
is a maximal order of $D$ contained in $\End_{\br{\FF}_3}(C)$.  
Since the order 
$\End(C)$ is
maximal, we see that $R = \End_{\br{\FF}_3}(C)$.  
Since the generators of $R$ are
defined over $\FF_9$, we see that all of the endomorphisms of $C$ are
defined over $\FF_9$.
\qed \end{pf*}

The endomorphism ring $\End(C^\wedge)$ of the formal group $C^{\wedge}$ is the
$p$-completion of $\End(C^\wedge)$.  This follows from a theorem of Tate
(see \cite{WaterhouseMilne}), which states that the $p$-completion
of the endomorphism ring of an
abelian variety is the the endomorphism ring of its $p$-divisible group.
However, in the case of the elliptic curve $C$, this may be deduced
in a more explicit manner.

Recall that the endomorphism ring of a height $2$ formal group over
$\FF_{p^2}$ is the
unique maximal order $\mathcal{O}_p$ of the $\QQ_p$-division algebra $D_p$ 
of invariant $1/2$
\cite[Appendix~B]{Ravenel}.  The ring $\mathcal{O}_p$ admits the
presentation
$$ \calO_p = \WW \langle S \rangle/ (S^2 = p, w^\sigma S = S w ) $$
where $\WW = \WW(\FF_{p^2})$ is the Witt ring of $\FF_{p^2}$, 
$w$ ranges over the
elements of $\WW$, and $\sigma$ is the lift of the
Frobenius $\Frob_p$ on $\FF_{p^2}$.

Since $C$ is supersingular, the formal group $C^\wedge$ is a height $2$
formal group.
The formal
group law $C^\wedge$ is seen, using the formulas of \cite{Silverman}, to have
$3$-series
$$ [3]_{C^\wedge}(T) = -T^9. $$

The Witt ring $\WW$ is the ring of integers of the unramified
quadratic extension of $\QQ_3$.
We may therefore identify the Witt ring $\WW$ with the subring $\ZZ_3[t]$ of
$\calO \otimes \ZZ_3$.  Let $\omega \in \FF_9$ be a primitive $8$th root of
unity.  We shall also let $\omega$ represent the Teichm\"uller lift to
$\WW$.  
Relation~\ref{eq:phitrel} implies that in $\calO \otimes \ZZ_3$ we have
\begin{equation}\label{eq:phiomegarel}
F \omega = \omega^3 F.
\end{equation}
Define $S$ to be the element $\omega^{-1}F \in \calO \otimes
\ZZ_3$.  Relations~\ref{eq:phiomegarel} and \ref{eq:phirel} imply that
we have
\begin{align}
S^2 & = 3 \label{eq:Srel} \\
S \omega & = \omega^3 S. \label{eq:Somegarel}
\end{align}
The following proposition follows immediately.

\begin{myprop}
There is an isomorphism of $\WW$-algebras
$$ \calO \otimes \ZZ_3 = \End(C) \otimes \ZZ_3 = \End(C^\wedge) \rightarrow 
\calO_3 $$
given by sending $F$ to $\omega S$.
\end{myprop}

The Morava stabilizer group $\MS = \MS_2$ is the group of units of $\calO_3$
$$ \aut(C^\wedge) = \calO_3^\times. $$

To understand the extended automorphism groups
\begin{align*}
G_{24} = \aut_{/\FF_3}(C) & = \aut(C) \rtimes Gal \\
\GG = \aut_{/\FF_3}(C^\wedge) & = \aut(C^\wedge) \rtimes Gal \\
\end{align*}
we simply must understand
the action of the Galois group $Gal = Gal(\FF_9/\FF_3)$.
Let $\sigma = \Frob_3$ denote the generator of $Gal$.
The Galois action on
$\End(C)$ is governed by the equations
\begin{align*}
\sigma^* F & = F \\
\sigma^* t & = -t.
\end{align*}
We deduce that the Galois action on $\aut(C)$ is given by
\begin{align*}
\sigma^* s & = s \\
\sigma^* t & = -t.
\end{align*}
The Galois action on $\WW \subset \aut(C^\wedge)$ is given by the 
lift of $\sigma$.
The Galois action on $\aut(C^\wedge)$
is then determined by the action of $\sigma$ on $S$.  We have
$$ \sigma^* S = \omega^{-2} S. $$

\begin{myrem}\label{rmk:Galoisaction}
This formula for the Galois action on $\MS$ differs from that that is
in more common usage (see, for instance, \cite{DevinatzHopkinsaction}).
The reason for this is that the Morava stabilizer group $\MS$ 
is usually taken to
be the group of automorphisms of the Honda height $2$ formal group $F_2$ over
$\FF_{p^2}$.  Since the formal groups $F_2$ and $C^\wedge$ are isomorphic
over $\br{\FF}_3$, they have isomorphic automorphism groups.  However,
because they are not isomorphic over $\FF_9$, their automorphism groups 
have non-isomorphic Galois actions.
\end{myrem}

\subsection{Effect on Morava modules}\label{sec:Morava}

In this section we summarize the effects of our operations on
$E_2$-homology.  
This is an adaptation of the techniques of Mahowald and Rezk
\cite{MahowaldRezk}.
Throughout this section, we freely use the notation established in
Section~\ref{sec:aut}. We shall denote the Morava
$E$-theory spectrum $E_2$ by $E$.  We are taking the spectrum $E$ to be 
the Goerss-Hopkins-Miller
$E_\infty$ ring spectrum associated to the formal group 
$C^\wedge$ over $\FF_9$.
When we write the $E$-homology group $E_*X$, 
we mean the Morava module $\pi_*(L_{K(2)}(E \wedge X))$.

Recall \cite{DevinatzHopkins} that for a closed subgroup
$F \subseteq \GG$ we have an isomorphism of Morava modules
$$ E_*E^{hF} \cong \map^c(\GG/F, E_*). $$

Let $C$ be the supersingular elliptic curve of Section~\ref{sec:aut}.
There is a
lift of this curve to a curve $\td{C}$ over $E^0 = \WW [[u_1]]$.  
$$ \td{C} : y^2 = 4x^3 + u_1 x^2 + 2x. $$

\begin{myprop}
The formal group law $\td{C}^\wedge$ is a universal deformation of
$C^\wedge$.
\end{myprop}

\begin{pf*}{Proof.}
We just need to see that the universal map classifying our deformation
induces an isomorphism on the base ring $E^0$.  This follows from the fact
that the $3$-series of $\td{C}^\wedge$ is given by \cite{Rezk}
$$ [3]_{\td{C}^\wedge}(T) \equiv u_1 T^3 + 2(1 + u_1^2)T^9 + \cdots \pmod
3. $$
\qed \end{pf*}

\begin{myrem}
The Serre-Tate theorem gives an equivalence between the category of
deformations of $C$ and the category of deformations of the formal group
$C^\wedge$.  Because $\td{C}^\wedge$ is a universal deformation of
$C^\wedge$, we may deduce that $\td{C}$ is a universal
deformation of $C$.
\end{myrem}

We define finite subgroups of $\GG$ to be the subgroups generated by the
following elements.
\begin{align*}
G_{24} & = \bra{s,t,\sigma} \\
SD_{16} & = \bra{\omega, \sigma} \\
D_8 & = \bra{t,\sigma}
\end{align*}
Work of Hopkins-Miller \cite{RezkHopkinsMiller} and Goerss-Hopkins
\cite{GoerssHopkins} shows that the $\GG$ action on $E_*$ lifts to an
$E_\infty$ action on the spectrum $E$.
We have the following spectra as homotopy fixed point
spectra.
\begin{align*}
E^{hG_{24}} & = TMF \\
E^{hSD_{16}} & = E(2) \\
E^{hD_8} & = TMF_0(2)
\end{align*}
The relationship to some of these groups and the curve $C$ is given below.
\begin{align*}
\GG & = \aut_{/\FF_3}(C^\wedge) \\
G_{24} & = \aut_{/\FF_3}(C) \\
D_8 & = \aut_{/\FF_3}(C,H).
\end{align*}
The subgroup 
$H$ is the level $2$ structure generated by the point $(0,0)$ of $C$.

\begin{myrem}\label{rmk:K(2)localTMF}
We pause to comment on the relationship between $TMF$ and $EO_2$.
The pair $(E,\td{C})$ is an $E_\infty$ elliptic spectrum \cite{HopkinsICM}.
This $E_\infty$ elliptic spectrum structure 
is classified by an $E_\infty$ ring map
$$ \kappa : TMF \rightarrow E. $$
Because $\td{C}$ is a universal deformation, 
the action of $G_{24}$ on $C$ extends to an action of $G_{24}$ on
$\td{C}$ (covering the action of $G_{24}$ on $E_*$).
The classifying map $\kappa$ therefore lifts to an $E_\infty$ ring map
$$ TMF \rightarrow E^{hG_{24}} = EO_2. $$
This map is a $K(2)$-local equivalence.  The essential point is that the 
$K(2)$-localization of $TMF$ will be given as the sections of the sheaf 
$\calO_{ell}$ over a formal neighborhood of the unique supersingular point
$C$ of the ($3$-local) moduli stack $\calM$ 
(see Section~\ref{sec:sheafrealize}).
\end{myrem}

Given an endomorphism
$$ \phi: C \rightarrow C $$
we get an induced endomorphism
$$ \widehat{\phi} : C^\wedge \rightarrow C^\wedge. $$
Assuming that $\widehat{\phi}$ is an automorphism, 
we may regard $\widehat{\phi}$ as an element of $\GG$, thus giving an 
$E_\infty$ ring map
$$ \widehat{\phi}: E \rightarrow E. $$
The induced map $E_* \widehat{\phi}$ on $E$-homology is given by 
$$
\xymatrix{
E_*E \ar[r]^{E_* \widehat{\phi}} \ar@{=}[d] &
E_*E \ar@{=}[d]
\\
\map^c(\GG,E_*) \ar[r]_{R_{\widehat{\phi}}^*} & 
\map^c(\GG,E_*)
}
$$
where $R_{\widehat{\phi}}$ is right multiplication by $\widehat{\phi}$.  If
$\widehat{\phi}$ is in the normalizer $N_\GG F$ of a closed subgroup $F$, then
there is an induced map
$$ \widehat{\phi} : E^{hF} \rightarrow E^{hF}. $$
Being in the normalizer implies that \emph{right} multiplication descends
to the right coset spaces $\GG/F$.  Thus the induced map on $E$-homology is
given by
$$ \map^c(\GG/F,E_*) \xrightarrow{R_{\widehat{\phi}}^*} \map^c(\GG/F,E_*).
$$

We may apply this discussion to translate our maps $\phi_f$, $\psi_d$, and
$\psi_{[2]}$
into the language of homotopy fixed point spectra.
The map $\phi_f^*$ is clearly the map induced from inclusion
$\iota_{D_8} : D_8 \hookrightarrow G_{24}$ giving
$$ \phi_f^* = Res_{D_8}^{G_{24}} : E^{hG_{24}} \rightarrow E^{hD_8}. $$
Consider the quotient isogeny
$$ \phi_H : C \rightarrow C/H $$
where $H$ is the level $2$ structure on $C$ generated by the point 
$(x,y) = (0,0)$.  The computations of Section~\ref{sec:psi_d} indicate that
the curve $C/H$ is given by the Weierstrass equation
$$ C/H : y^2 = x^3 + x. $$
Note that this differs from our Weierstrass equation for $C$, but since
there is only one isomorphism class of supersingular curve at $p=3$, this
curve must be isomorphic to $C$.  Indeed, we find it convenient to take the
isomorphism
$$ \mu_{\omega^{-1}} : C/H \rightarrow C $$
given by
\begin{align*}
x & \mapsto \omega^{-2} x \\
y & \mapsto \omega^{-3} y. 
\end{align*}
Thus we may consider the composite 
$$ \psi_H : (C,H) \xrightarrow{\phi_H} C/H \xrightarrow{\mu_{\omega^{-1}}}
C. $$
The induced map $\widehat{\psi}_H$ on the formal group is readily
computed from the formulas of Section~\ref{sec:phi_q}.  One finds that,
regarding $\widehat{\psi}_H$ as an element of $\GG$, we have
$$ \widehat{\psi}_H = 1+t. $$
The group element $1+t$ is in the normalizer $N_\GG D_8$.
Thus the map $\psi_d^*: TMF(2) \rightarrow TMF(2)$ is given by
$$ \psi_d^* = [1+t] : E^{hD_8} \rightarrow
E^{hD_8}. $$

Since $\psi_{[2]}$ corresponds to the second power isogeny, it is given by
the element $2 \in \GG$, so that
$$ \psi_{[2]} = [2] : E^{hG_{24}} \rightarrow E^{hG_{24}}. $$

We could have used $1-t \in \GG$ to represent the map $\psi_d^*$.  
It is irrelevant on
$E^{hD_8}$ since the cosets $(t+1)D_8$ and $(1-t)D_8$ are equal in
$\GG/D_8$. The formula
$$ (t+1)(1-t) = 2 $$
in $\GG$ implies relation~\ref{eq:rel2}
$$ (\psi_d^*)^2 = \psi_{[2]} : E^{hD_8} \rightarrow E^{hD_8}$$
required in the construction of $Q(2)$ given in
Section~\ref{sec:superrealize}.

\subsection{Calculation of $V(1)_*(Q(2))$ at $p=3$}
\label{sec:Q(2)V(1)}

Let $V(1)$ be the $3$-primary Smith-Toda complex.
In this section we will compute the $V(1)$-homology group 
$V(1)_*(Q(2))$.  We shall compare this computation 
to $V(1)_*(S)$, where $S$ is the
($K(2)$-local) sphere.

\subsubsection{Computation of the ANSS}

We shall compute the Adams-Novikov $E_2$ term for $Q(2) \wedge V(1)$ using
the hypercohomology 
spectral sequence~\ref{eq:ANE2SS} (modulo the ideal $(3,v_1)$).  
We first describe the $E_1$ term.

We recall \cite{Rezk} that the the $[3]$-series of the formal group law of the 
elliptic curve $C_\mathbf{q}$
has
$$ v_1 \equiv q_2 \pmod 3. $$
Therefore, the ANSS for $TMF \wedge V(1)$ is the cohomology of the quotient
Hopf algebroid of $(\br{B},\Gamma_{\br{B}})$ (see Section~\ref{sec:ANcomplex}) 
given by 
$$ (\br{B}/I_2, \Gamma_{\br{B}}/I_2 ) $$
where $I_2$ is the invariant regular ideal $(3,q_2) \subset \br{B}$.
The ANSS $E_2$ term for $TMF \wedge V(1)$ is computed \cite{Bauer}, \cite{Rezk}
to be
$$ E_2(TMF\wedge V(1)) = 
H^*(\Gamma_{\br{B}}/I_2) = \FF_3[q_4^{\pm 1}, \beta] \otimes E[\alpha] $$
where the generators have bidegrees
\begin{align*}
\abs{q_4} & = (0,8) \\
\abs{\beta} & = (2,12) \\
\abs{\alpha} & = (1,4). 
\end{align*}
We have $v_2 = -q_4^2$.  
The ANSS $E_2$ term for $TMF_0(2) \wedge V(1)$ is concentrated on the zero
line, and is given by
$$ E_2(TMF_0(2) \wedge V(1)) = \FF_3[q_4^{\pm 1}]. $$  

The differentials in the complex $C^*(Q(2) \wedge V(1))$ of
Section~\ref{sec:ANcomplex} are computed from the formulas of
Section~\ref{sec:mapscomp}.  Namely, we have
\begin{align*}
\phi_f^*(q_4^k) & \equiv (q_4)^k \pmod {I_2} \\
\psi_{[2]}^* (q_4^k) & \equiv q_4^k \pmod {I_2} \\
\psi_d^* (q_4^k) & \equiv (-q_4)^k \pmod {I_2} \\
\phi_q^*(q_4^k) & \equiv (-q_4)^k \pmod {I_2}.
\end{align*}
The
non-zero differentials $D_i$ in the total complex $C^*(Q(2)\wedge V(1))$ 
are given by the following formulas.
\begin{align*}
D_0 : C^0(\Gamma_{\br{B}}) & \rightarrow C^1(\Gamma_{\br{B}}) \oplus
C^0(\Gamma_{\br{B}}) \oplus MF_0(2)
\\
\\
D_0(q_4^k) & \equiv 
(0,0,q_4^k) \: \pmod {I_2}, \quad k \equiv 1 \pmod 2
\\
\\
D_1: C^1(\Gamma_{\br{B}}) \oplus
\br{C}^0(\Gamma_{\br{B}}) \oplus MF_0(2) & \rightarrow 
C^2(\Gamma_{\br{B}}) \oplus
\br{C}^1(\Gamma_{\br{B}}) \oplus MF_0(2)
\\
\\
D_1(0,0,q_4^k) & \equiv (0,0,-q_4^k) \pmod {I_2}, \quad k \equiv 0 \pmod 2
\\
D_1(0,q_4^k,0) & \equiv (0,0,-q_4^k) \pmod {I_2}, \quad k \in \ZZ
\end{align*}

We therefore have the following proposition.

\begin{myprop}
The $E_2$-term of the ANSS for $Q(2) \wedge V(1)$ is given by
$$ \FF_3[\beta, v_2^{\pm 1}]\otimes E[\zeta] \{1,\alpha, h_1, b_1 \}. $$
\end{myprop}

The generators are given in the following table.  The representative
is given the name of the corresponding
element of the ANSS $E_2$-term of the layer of the tower 
(The brackets refer
to the level of the tower that the generator lives in). 
\vspace{.2in}

\begin{center}
\begin{tabular}{c|c|c}
\hline
Generator & Representative & Bidegree \\
\hline \hline
$\beta$ & $\beta[0]$ & $(2,12)$ \\
$v_2$ & $-q_4^2[0]$ & $(0,16)$ \\ 
$\alpha$ & $\alpha[0]$ & $(1,4)$ \\
$h_1$ & $\alpha q_4[0]$ & $(1,12)$ \\
$b_1$ & $\beta q_4^3[0]$ & $(2,36)$ \\
$\zeta$ & $1_{MF}[1] - 1_{MF_0(2)}[1]$ & $(1,0)$ \\
\hline
\end{tabular}
\end{center}
\vspace{.2in}

The differentials in the ANSS for $Q(2)\wedge V(1)$ 
follow from the differentials in
the ANSS for $TMF\wedge V(1)$ (see \cite{Bauer}, \cite{Rezk}).

\begin{myprop}\label{prop:Q(2)V(1)ANSSdiffs}
The differentials in the ANSS for $Q(2) \wedge V(1)$ are given by
\begin{align*}
d_5(v_2^k) & = v_2^{k-2} \beta^2 h_1 , \quad k \equiv 2,3,4,6,7,8 \pmod 9 \\
d_9(v_2^k \alpha) & = v_2^{k-3} \beta^5 , \quad k \equiv 3,4,8 \pmod 9 \\
d_9(v_2^k h_1) & = \beta^4 v_2^{k-4}, \quad k \equiv 3,7,8 \pmod 9 \\
d_5(v_2^k b_1) & = v_2^k \beta^3 \alpha, \quad k \equiv 0,1,2,5,6,7 \pmod 9
\end{align*}
and these differentials are propagated freely by multiplication by $\beta$
and $\zeta$.  This completely describes the ANSS.
\end{myprop}

\subsubsection{The $V(1)$-homology groups $V(1)_*(Q(2))$ and
$V(1)_*(S)$}\label{sec:piQ(2)V(1)}

\begin{figure}
\includegraphics[width=4.9in]{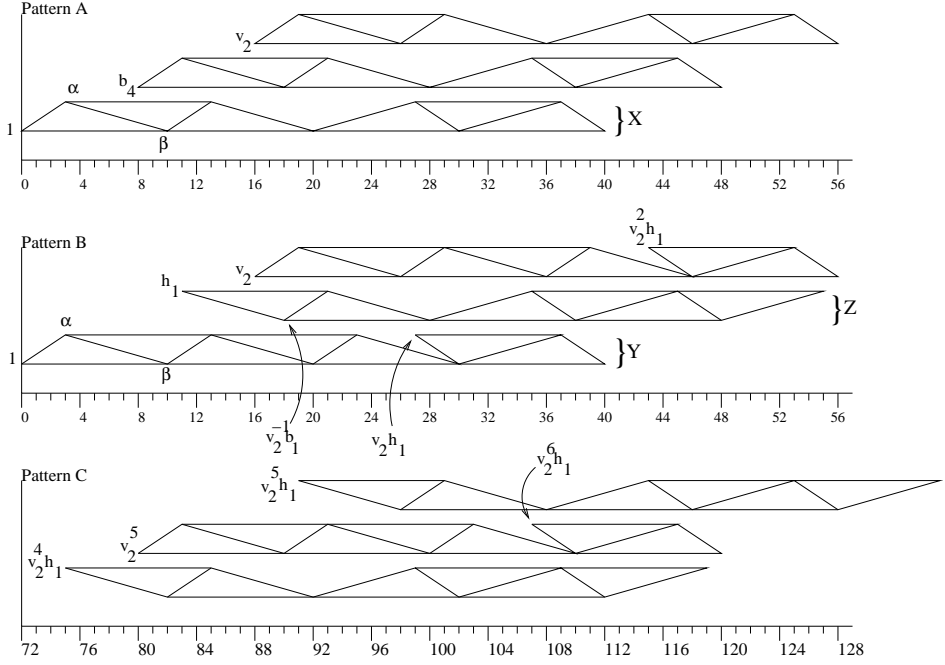}
\caption{Graphical depiction of patterns $A$, $B$, and $C$.}
\label{fig:patterns}
\end{figure}

Define patterns $A$, $B$, and $C$ of homotopy groups as follows
\begin{align*}
A & = X \otimes \FF_3\{1, b_4, v_2 \} \\
B & = Y \otimes \FF_3\{1, v_2 \} \oplus Z \\
C & = Z \otimes \FF_3\{v_2^4, v_2^5 \} \oplus Y \otimes \FF_3\{v_2^5 \}
\end{align*}
where $X$, $Y$, and $Z$ are the following graded $\FF_3$-vector spaces.
\begin{align*}
X & = \FF_3 \{ 1, \alpha, \beta, \alpha \beta, \beta^2, \bra{\alpha, \alpha,
\beta^2}, \beta^3, \bra{\alpha, \alpha, \beta^3}, \beta^4 \} \\
Y & = \FF_3 \{1, \alpha, \beta, \alpha\beta, \beta^2, \alpha\beta^2,
v_2h_1, \beta^3, \beta v_2h_1, \beta^4 \} \\
Z & = \FF_3 \{ h_1, v_2^{-1}b_1, \beta h_1, \beta v_2^{-1}b_1, \bra{\alpha,
\alpha, \beta v_2^{-1}b_1}, \beta^2 v_2^{-1}b_1, \bra{\alpha, \alpha,
\beta^2 v_2^{-1}b_1}, \\
& \qquad \beta^3 v_2^{-1}b_1, \bra{\alpha, \alpha, \beta^3 v_2^{-1}b_1} \}
\end{align*}
Figure~\ref{fig:patterns} gives a graphical description of the patterns
$A$, $B$, and $C$, and the patterns $X$, $Y$, and $Z$ sit inside of them as
indicated.  In this figure, the lines of the appropriate length 
depict multiplication
by $\alpha$, multiplication by $\beta$, and application of the Toda
bracket $\langle \alpha, \alpha, - \rangle$.

Let $Q'(2)$ be the fiber of the unit
$$ Q'(2) \rightarrow S \xrightarrow{\eta} Q(2). $$
With these patterns, we have the following computations.  The first two 
computations are given in \cite{GoerssHennMahowald}.  The 
computation of $V(1)_*(Q(2))$ is an
immediate consequence of Proposition~\ref{prop:Q(2)V(1)ANSSdiffs}.
\begin{gather*}
V(1)_*(TMF) = (A \oplus \Sigma^{72} A)\otimes P[v_2^{\pm 9}] 
\\
V(1)_*(S) = ((B \oplus C) \oplus \Sigma^{29} (B
\oplus C)^{\vee})\otimes P[v_2^{\pm 9}] \otimes E[\zeta] 
\\
V(1)_*(Q(2)) = (B \oplus C)\otimes P[v_2^{\pm 9}]
\otimes E[\zeta] 
\\
V(1)_*( Q'(2)) = \Sigma^{29} (B
\oplus C)^{\vee}\otimes P[v_2^{\pm 9}] \otimes E[\zeta]
\end{gather*}
The computation of the $V(1)$-homology of $Q'(2)$ follows from the fact
that the computations of the $V(1)$-homology of $S$ and $Q(2)$
imply there is a short exact sequence
$$ 0 \rightarrow V(1)_*(Q'(2)) \rightarrow V(1)_*(S) 
\rightarrow V(1)_*(Q(2)) \rightarrow 0 $$

Thus it appears that $V(1)$ is built out of
$Q(2) \wedge V(1)$ and a $28^\mathrm{th}$ 
suspension of the Brown-Comenetz dual of $Q(2) \wedge V(1)$.  This is the key
computation to our proof that $Q'(2) \simeq D(Q(2))$
(Theorem~\ref{thm:cofiber}).  Gross-Hopkins duality \cite{HopkinsGross} 
implies that after
smashing with $V(1)$, Spanier-Whitehead duality and Brown-Comenetz duality
agree up to a suitable suspension in the $K(2)$-local category.

\section*{Part 2: The $3$-primary $K(2)$-local sphere}
\addcontentsline{toc}{part}{Part 2: The $3$-primary $K(2)$-local sphere
\hfill}
\addtocounter{section}{1}
\setcounter{subsection}{0}
\setcounter{mythm}{0}

Let $K$ denote $K(2)$, $E$ denote $E_2$, and $Q$ denote
$Q(2)$.  The ($K$-local) Spanier-Whitehead dual of $X$ will be denoted by 
$$ DX = F(X, L_K S). $$ 
Recall that the $n$th monochromatic layer $M_nX$ is the fiber
$$ M_nX \rightarrow L_nX \rightarrow L_{n-1}X. $$
The Gross-Hopkins dual $I_nX$ is the Brown-Comenetz dual of $M_nX$.
We
shall let $I_KX$ denote the
$K$-localization of $I_nX$.  In particular, we shall write $I_K$ for $I_KS$.  
The spectrum $I_K$ is
invertible \cite{HopkinsGross}, \cite{HoveyStrickland}.
Therefore, we have \cite{Strickland}
\begin{equation}\label{eq:dualizing}
I_KX \simeq DX \wedge I_K.
\end{equation}
In this part of the paper we wish to prove the following theorem.

\begin{mythm}\label{thm:cofiber}
There is a cofiber sequence
$$ DQ \xrightarrow{D \eta} S \xrightarrow{\eta} Q $$
where $\eta$ is the unit of the ring spectrum $Q$.
\end{mythm}

By equation~\ref{eq:dualizing}, we have $DQ \simeq I_KQ \wedge
I_K^{-1}$.  We are able to use this equation to give an explicit
description of $DQ$.  We shall prove

\begin{mythm}\label{thm:DQ}
The spectrum $DQ$ refines a tower of spectra of the form
$$ \Sigma^{46}TMF_0(2) \rightarrow \Sigma^{46}(TMF \vee TMF_0(2)) 
\rightarrow \Sigma^{46}TMF. $$
\end{mythm}

It will turn out that the spectra $TMF$ and $TMF_0(2)$ are Gross-Hopkins
self-dual, up to a suitable suspension.  
This will follow from the relationship of
the Gross-Hopkins dual to the Mahowald-Rezk dual.  The Mahowald-Rezk duals
of certain $fp$-spectra are easily computed
using the theory of \cite{MahowaldRezk}.  The number $46$ in
Theorem~\ref{thm:DQ} arises from $46 = 22+22+2$ where the first two $22$'s
arise from the equations
\begin{gather*}
I_KTMF \simeq \Sigma^{22} TMF \quad \text{(Proposition \ref{prop:I_KTMF})}
\\
I_K^{-1} \wedge TMF \simeq \Sigma^{22} TMF \quad \text{(Proposition
\ref{prop:I_K^TMF})}
\end{gather*}
and the final factor of $2$ arises from the fact that the tower 
$$ TMF \rightarrow TMF \vee TMF_0(2) \rightarrow TMF_0(2) $$
has length $2$ (Lemma~\ref{lem:shiftdual}).  
Thus we recover, at least in form, the tower of
\cite{GHMR}.  The maps in the tower for $DQ$ 
given in Theorem~\ref{thm:DQ} are the
Spanier-Whitehead duals of the maps in the tower for $Q$.

We shall first give a general abstract framework for duality pairings 
in triangulated closed symmetric monoidal categories, where we define the
notion of a Lagrangian for a hyperbolic pairing.
We shall prove that for the duality pairing for 
the $K(2)$-local sphere is hyperbolic, in
the sense that $S$ is built from two dual spectra.  
The Lagrangian for this decomposition will be $\br{S} = E^{h\GG^1}$, where
$\GG^1$ is the kernel of the reduced norm $\GG^1 \rightarrow \ZZ_3$.
The decomposition we consider is well known \cite{DevinatzHopkins}, but
this particular interpretation of it is new.  
We then use this decomposition to
prove that $Q$ is also a Lagrangian for $S$, thus proving
that $S$ may also be built out of $Q$ and $DQ$ (Theorem~\ref{thm:cofiber}).

The second part of this paper is organized as follows.  
In Section~\ref{sec:pairings}, we give our abstract duality framework.
In
Section~\ref{sec:Bousfield}, it is recalled that $V(1)$-homology
equivalences are $K(2)$-local equivalences.  In Section~\ref{sec:Sbar} we
prove Proposition~\ref{prop:Sbar}, which says  
that for arbitrary $n$, 
the canonical pairing on $L_{K(n)}S$ is hyperbolic with 
Lagrangian $\br{S}$.
Section~\ref{sec:I_KTMF} identifies $I_KTMF$ and Section~\ref{sec:I_K^TMF}
identifies $I_K \wedge TMF$.  In Section~\ref{sec:DQ} we use these
identifications to prove Theorem~\ref{thm:DQ}.  In Section~\ref{sec:Qbar}
we define a spectrum $\br{Q}$, and in Section~\ref{sec:Sbarlagrangian} we
prove that $\br{Q}$ is a Lagrangian for $\br{S}$.  This uses some
decompositions of mapping spaces of homotopy fixed point spectra that are
recalled from \cite{GHMR} in Section~\ref{sec:GHMR}.
Section~\ref{sec:Qdecomp} describes how $Q$ may be built from two copies of
$\br{Q}$.  Theorem~\ref{thm:cofiber} is finally proved in
Section~\ref{sec:cofiberproof}.

\subsection{Duality pairings}\label{sec:pairings}

In this section we briefly outline the abstract framework in which the
duality phenomena of this paper reside.  Compare with \cite{MayPicard},
\cite{FauskHuMay}.

We first provide some motivation.  Suppose that $V$ is a finite dimensional
vector space over a field $k$. A (bilinear) pairing is a linear
homomorphism
$$ \alpha: V \otimes V \rightarrow k. $$
It is a perfect pairing if the adjoint homomorphism
$$ \td{\alpha} : V \rightarrow V^* $$
is an isomorphism.  A Lagrangian for such a pairing is a subspace $W
\subseteq V$ such that the induced sequence
$$ W \xrightarrow{\iota} V \cong V^* \xrightarrow{\iota^*} W^* $$
is a short exact sequence, where the isomorphism between $V$ and $V^*$ is
given by $\td{\alpha}$.  A pairing for which a Lagrangian exists is
a hyperbolic pairing.  We wish to provide a framework of hyperbolic
pairings in the setting of triangulated closed symmetric monoidal categories.

Let $(\calC, \wedge, F(-,-), S)$ 
be a closed symmetric
monoidal category with compatible triangulated structure
\cite{May}.  
Let $DX$ denote the dual $F(X,S)$ of $X$.  One says that $X$ is 
\emph{reflexive} if the natural map
$$ \td{\epsilon}: X \rightarrow DDX $$
is an isomorphism.  This should be contrasted with the stronger 
notion of $X$ being
\emph{dualizable}, where one insists that the natural map
$$ DX \wedge X \rightarrow F(X,X) $$
be an isomorphism.

We shall refer to any map
$$ \alpha: X \wedge X \rightarrow I $$
where $I$ is an element of $Pic(\calC)$,
as a \emph{pairing} on $X$.  A pairing $\alpha$ will be called
\emph{perfect} if the adjoint map
$$ \td{\alpha}: X \rightarrow F(X,I) \cong DX \wedge I $$
is an isomorphism.  Note that it is immediate that 
if $X$ possesses a perfect pairing $\alpha$, then $X$ is reflexive.

Suppose that $X$ has a perfect pairing $\alpha: X\wedge X \rightarrow I$.
Given a map $l: Y \rightarrow X$, we may define a dual map $l^\vee$ 
by the composite
$$ l^\vee: X \xrightarrow[\cong]{\td{\alpha}} DX\wedge I
\xrightarrow{Dl\wedge 1} DY\wedge I. $$
We shall say that $l: Y \rightarrow X$ is a \emph{left Lagrangian} for
$\alpha$ if the sequence
$$ Y \xrightarrow{l} X \xrightarrow{l^\vee} DY\wedge I $$
extends to a distinguished triangle in $\calC$.
Alternatively, given a map $r: X \rightarrow Z$, we may define a map
$r^\vee$ to be the composite
$$ r^\vee : DZ\wedge I \xrightarrow{Dr \wedge 1} DX\wedge I
\xrightarrow[\cong]{\td{\alpha}^{-1}} X $$
We shall say that $r: X \rightarrow Z$ is a \emph{right Lagrangian} for
$\alpha$ if the sequence
$$ DZ\wedge I \xrightarrow{r^\vee} X \xrightarrow{r} Z $$
extends to a distinguished triangle of $\calC$.  We have the following
proposition which enumerates some useful facts regarding these definitions.

\begin{myprop}
Suppose that $X$ has a perfect pairing $\alpha: X\wedge X \rightarrow I$.
\begin{enumerate}

\item If $Y$ is reflexive, then given $l:Y \rightarrow X$, we have 
$l^{\vee \vee} \cong l$ under the isomorphism $D(DY \wedge I)\wedge I \cong
Y$.

\item If $Z$ is reflexive, then given $r:X \rightarrow Z$, we have 
$r^{\vee \vee} \cong r$ under the isomorphism $D(DZ \wedge I)\wedge I \cong
Z$.

\item If $l: Y \rightarrow X$ is a left Lagrangian, then $Y$ is reflexive, 
and $l^\vee: X \rightarrow DY\wedge I$ is a right Lagrangian.

\item If $r: X \rightarrow Z$ is a right Lagrangian, then $Z$ is reflexive, 
and $r^\vee: DZ\wedge I \rightarrow X$ is a left Lagrangian.

\end{enumerate}
\end{myprop}

\begin{pf*}{Proof.}
We shall only prove $(3)$.  The other parts use similar methods.  By
hypothesis, there exists a map $\delta$ so that the top row of the
following diagram is a distinguished triangle.
$$
\xymatrix@C+1em{
Y \ar[r]^l \ar[d]^{\epsilon_1} \ar@/_4pc/[dd]_{\td{\epsilon}} &
X \ar[d]_{\cong}^{\td{\alpha}} \ar[r]^{l^\vee} &
DY \wedge I \ar@{=}[d] \ar[r]^\delta &
\Sigma Y \ar[d]^{\Sigma \epsilon_1} 
\\
D(DY \wedge I)\wedge I \ar[r]_-{Dl^\vee \wedge 1} \ar[d]^\cong &
DX\wedge I \ar[r]_{Dl\wedge 1} &
DY \wedge I \ar[r]_-{D\Sigma^{-1}\delta\wedge 1} &
\Sigma D(DY\wedge I)\wedge I
\\
DDY
}
$$
The
second row is distinguished, since it is what one gets when one applies the
functor
$D(-)\wedge I$ to the first row, after we have turned the triangle once.
The sign on the map $D\Sigma^{-1}\delta\wedge 1$ is correct: we get one
factor of $(-1)$ from turning the triangle and one factor of $(-1)$ from
axiom TC2 of \cite{May}.
One concludes that $\epsilon_1$ is an
isomorphism, hence $\td{\epsilon}$ is.
\qed \end{pf*}

We conclude with a degenerate example.

{\bf Example.}
Given an element $I$ of $Pic(\calC)$, there is a unique pairing
$$ \ast \wedge \ast \rightarrow I $$
which is perfect. 
Suppose that $\ast \rightarrow Z$ is a right Lagrangian.  Then there exists
a map $\delta$ so that the
sequence
$$ DZ \wedge I \rightarrow \ast \rightarrow Z \xrightarrow{\delta}
\Sigma DZ \wedge I$$
is a distinguished triangle, which means that there is an isomorphism
$$ \td{\alpha} = \delta: Z \rightarrow DZ \wedge \Sigma I. $$
The adjoint
$$ \alpha: Z \wedge Z \rightarrow \Sigma I $$
is a perfect pairing on $Z$.  Conversely, given $Z$ with a perfect pairing
taking values in $\Sigma I$, we see that the unique map $\ast \rightarrow
Z$ is a right Lagrangian.  
The choice of perfect pairing on $Z$ coincides with the choice of how to
complete the sequence $DZ\wedge I \rightarrow \ast \rightarrow Z$ to a
distinguished triangle.
Indeed, given that $Z$ admits a perfect pairing,
the collection of perfect pairings on
$Z$ is (non-canonically) isomorphic to $\aut(Z)$, which is also the
automorphism group of the triangle
$$ \Sigma^{-1}Z \rightarrow \ast \rightarrow Z \xrightarrow{=} Z. $$

\subsection{$K(2)$-local equivalences}\label{sec:Bousfield}

In what follows we will often want to deduce certain maps are $K(2)$-local 
equivalences.  Of course, these are precisely maps which are
$K(2)$-homology equivalences.  
However, it is sometimes useful to use homotopy group
calculations instead.  A $K(1)$ version of the following lemma was used to
identify the $J$ spectrum with $L_{K(1)}S$ in \cite{Bousfield}.

\begin{mylem}\label{lem:Bousfield}
Suppose that $f:X \rightarrow Y$ is a map for which the
induced map
$$ f\wedge V(1) : X\wedge V(1) \rightarrow Y \wedge V(1) $$
induces an isomorphism of homotopy groups.  
Then $f$ is a $K(2)$-equivalence.
\end{mylem}

\begin{pf*}{Proof.}
Since $f \wedge V(1)$ is an equivalence, $f \wedge V(1) \wedge E(2)$ is
an equivalence.  But $V(1) \wedge E(2)$ is equivalent to $K(2)$.  
Thus $f$ is a $K(2)$
equivalence.
\qed \end{pf*}

\subsection{$L_{K(n)}S$ is hyperbolic}\label{sec:Sbar}

In this section we work at an arbitrary prime $p$, and 
prove that the canonical duality pairing 
for the $K(n)$-local sphere is hyperbolic.  
For the purposes of this section everything is
implicitly $K(n)$-local.
We abbreviate $K$ for $K(n)$ and $E$ for $E_n$.
Let $\MS$ be the $n$th Morava stabilizer group.  
There is a reduced norm map
$$ \MS \rightarrow \ZZ_p $$
whose kernel we shall identify as as $\MS^1$.  Let $\GG^1$ denote the
Galois extension $\MS^1 \rtimes Gal(\FF_{p^n}/\FF_p)$.  Let $\br{S}$ be the
homotopy fixed point spectrum $E^{h\GG^1}$.  
The spectrum $\br{S}$ is an $E_\infty$ ring spectrum \cite{GoerssHopkins}.
The group $\ZZ_p$ acts on $\br{S}$ via $E_\infty$ maps. 
Let $\tau$ be a topological generator of
$\ZZ_p$.  For convenience of group-ring notation, we shall 
use multiplicative notation for the additive group
$\ZZ_p$.  In \cite{DevinatzHopkins} it is proved that there is a
homotopy fiber sequence
\begin{equation}\label{eq:sbarfiber}
\Sigma^{-1} \br{S} \xrightarrow{\delta} S \xrightarrow{r} \br{S} 
\xrightarrow{\tau-1} \br{S}.
\end{equation}

Define a pairing $\gamma$ by the composite  
\begin{equation}\label{eq:gammadef}
\gamma: 
\br{S} \wedge \br{S} \xrightarrow{\mu} \br{S} \xrightarrow{\delta}
S^1
\end{equation}
where $\mu$ is the product on $\br{S}$.  We shall refer to its adjoint as
$\td{\gamma}$.
$$ \td{\gamma}: \br{S} \rightarrow \Sigma D\br{S} $$

The following proposition simultaneously states that 
$\gamma$ is a perfect pairing, and that
the map $S \rightarrow \br{S}$ 
is a (right) Lagrangian for $S$.

\begin{myprop}\label{prop:Sbar}
The map $\td{\gamma}$ is an equivalence
which makes the left square of the following diagram commute.
$$
\xymatrix{
\Sigma^{-1} \br{S} \ar[r]^{\delta} \ar[d]^\simeq_{\td{\gamma}} &
S \ar[r]^r \ar@{=}[d] &
\br{S} \ar@{=}[d]
\\
D\br{S} \ar[r]_{Dr} &
S \ar[r]_r &
\br{S}}
$$
In particular, the bottom row is a cofiber sequence.
\end{myprop}

For a spectrum $X$ and a profinite set $K = \varprojlim_i K_i$, let
$X[[K]]$ be defined to be the homotopy inverse limit
$$ X[[K]] = \varprojlim_i X \wedge (K_i)_+. $$
Note that the natural map
$$ X \wedge S[[K]] \rightarrow X[[K]] $$
is not in general an equivalence.

We shall first need the following theorem, which is a proven in
\cite{BehrensDavis}.  

\begin{mythm}\label{thm:FEHEH}
Suppose that $H$ and $K$ are closed subgroups of $\GG$.
Then there is an
equivalence
$$ \alpha: (E[[\GG/H]])^{hK} \xrightarrow{\simeq} F(E^{hH}, E^{hK}). $$
\end{mythm}

\begin{mycor}\label{cor:FEHEH}
Suppose that $N$ is a closed normal subgroup of $\GG$.  Then there is an 
equivalence
$$ \alpha: E^{hN}[[\GG/N]] \xrightarrow{\simeq} F(E^{hN}, E^{hN}). $$
\end{mycor}

\begin{mycor}\label{cor:Gaction}
Suppose that $N$ is a closed normal subgroup of $\GG$.  Then there is an
action map
in the $K(n)$-local stable homotopy category
$$ \xi : S[[\GG/N]] \wedge E^{hN} \rightarrow E^{hN} $$
with the property that for every coset $gN \in \GG/N$, the composite
$$ \{ gN \}_+ \wedge E^{hN} \rightarrow S[[\GG]] \wedge E^{hN}
\xrightarrow{\xi} E^{hN} $$
coincides with the action of $gN$ on $E^{hN}$.
\end{mycor}

\begin{pf*}{Proof.}
Let $\td{\alpha}$ be the adjoint of the map $\alpha$ of
Corollary~\ref{cor:FEHEH}
$$ \td{\alpha}: E^{hN}[[\GG/N]] \wedge E^{hN} \rightarrow E^{hN}. $$
Then the map $\xi$ is given by the composite
$$ \xi : S[[\GG/N]] \wedge E^{hN} \xrightarrow{\eta \wedge 1 \wedge 1}
E^{hN}[[\GG/N]] \wedge E^{hN} \xrightarrow{\td{\alpha}} E^{hN} $$
where $\eta$ is the unit of the ring spectrum $E^{hN}$.
\qed \end{pf*}

\begin{mylem}\label{lem:finiteness}
For any closed subgroup $H$ of $\GG$, and profinite set $K = 
\varprojlim_i K_i$ for which all of the maps in the tower $\{ K_i \}$ are
surjections,
the natural map
$$ S[[K]] \wedge E^{hH} \rightarrow [[K]]E^{hH} $$
is an equivalence.
\end{mylem}

The proof of Lemma~\ref{lem:finiteness} will requires the following
finiteness result.

\begin{mylem}\label{lem:nilpotent}
There exists a dualizable spectrum $Y$ so that 
\begin{enumerate}
\item
$Y$-homology isomorphisms are equivalences.

\item
For every closed subgroup $H$ of $G$,
the $Y$-homology groups $Y_*E^{hH}$ are
finite in each degree.
\end{enumerate}
\end{mylem}

\begin{pf*}{Proof.}
Let $M$ be a type $n$ complex and let $U$ be an open normal 
subgroup of $\GG$ of
finite cohomological dimension (such a subgroup exists since $\GG$ is
compact $p$-adic analytic).
Let $Y$ be the spectrum $E^{hU} \wedge M$.  The spectrum $Y$ is dualizable
since the spectrum $E^{hU}$ is
dualizable.  In fact, the spectrum $E^{hU}$ is self-dual \cite{Rognes}.

{\it Proof of (1).} 
It suffices to show that for any spectrum $X$, if $Y \wedge X$ is
null, then $X$ must be null.  
In \cite{Rognes} it is shown that the norm map is an equivalence
$$ N: (E^{hU})_{h\GG/U} \xrightarrow{\simeq} (E^{hU})^{h\GG/U}. $$
Combined with the fact that $E^{hU}$ is self-dual,
we have
\begin{align*}
(Y \wedge X)^{h\GG/U} & \simeq F(E^{hU},X)^{h\GG/U} \wedge M \\
& \simeq F((E^{hU})_{h\GG/U},X) \wedge M \\
& \simeq F((E^{hU})^{h\GG/U},X) \wedge M \\
& \simeq F(E^{h\GG},X) \wedge M \\
& \simeq M \wedge X.
\end{align*}
Therefore $M \wedge X$ is null, but this implies that $X$ is null.

{\it Proof of (2).}
Using Theorem~\ref{thm:FEHEH}, and the fact that $E^{hU}$ is dualizable and
self-dual, we have
\begin{align*}
Y \wedge E^{hH} & = E^{hU} \wedge M \wedge E^{hH} \\
& \simeq F(E^{hU}, E^{hH})\wedge M \\
& \simeq (E[\GG/U])^{hH} \wedge M.
\end{align*}
As a left $H$-set, the finite $\GG$-set $\GG/U$ breaks up into a finite set of
$H$-orbits
$$ \GG/U \cong \coprod_{\GG/UH} H/(H \cap U). $$
Since $H/(H \cap U)$ is finite, there is an $H$-equivariant equivalence
$$ E[\GG/U] \simeq \bigvee_{\GG/UH} \map(H/(H \cap U), E). $$
By Shapiro's lemma (see, for example, \cite{BehrensDavis}), 
we deduce that the $H$-homotopy fixed points are
given by
$$ (E[\GG/U])^{hH} \simeq \bigvee_{\GG/UH} E^{hH\cap U}. $$
The homotopy fixed point spectral sequence gives a spectral sequence
$$ E_2^{s,t} = \bigoplus_{\GG/UH} H^s_c(H \cap U; E_tM) \Rightarrow
Y_{t-s}(E^{hH}). $$
The group $H \cap U$, being a subgroup of $U$, has finite cohomological
dimension \cite[4.1.2]{SymondsWeigel}, so the spectral sequence has a
horizontal vanishing line.  
Since $M$ is a type $n$ complex, the $E$-homology of $M$ is finite in each
degree.
Since the group $H \cap U$ is a
closed subgroup of the $p$-adic analytic group $\GG$, we may conclude that 
the cohomology 
groups $H^s_c(H \cap E_tM)$ are finite for each $s$ and $t$ (combine
Theorem~5.1.2, Lemma~5.1.4, and Proposition~4.2.2 of \cite{SymondsWeigel}).
We deduce that the $Y$-homology $Y_*E^{hH}$ is finite in each degree.
\qed \end{pf*}

\begin{pf*}{Proof of Lemma~\ref{lem:finiteness}.}
It suffices to show that the map
$$ S[[K]] \wedge E^{hH} \wedge Y \rightarrow 
[[K]]E^{hH} \wedge Y $$
is an equivalence, where $Y$ is the dualizable spectrum of
Lemma~\ref{lem:nilpotent}.  
The homotopy fixed point spectrum $E^{hH}$ is a colimit of homotopy fixed
point spectra 
$$ E^{hH} = \varinjlim_{H \le V \le_o \GG} E^{hV} $$
where the groups $V$ are open \cite{DevinatzHopkins}.
The spectra $E^{hV}$ are dualizable \cite{HoveyStrickland}.
Since homotopy inverse limits commute with smash products with dualizable
spectra, we are reduced to showing that the natural map
$$ \varinjlim_{H \le V \le_o \GG} \varprojlim_{i} ([K_i]E^{hV} \wedge Y) 
\rightarrow 
\varprojlim_{i} \varinjlim_{H \le V \le_o \GG} ([K_i]E^{hV}\wedge Y)
$$
is an equivalence.  By Theorem~3.7 of \cite{Mitchell}, it suffices to show
that the induced maps
\begin{equation}\label{eq:limcolim}
\varinjlim_{H \le V \le_o \GG} {\varprojlim_i}^s \pi_*(E^{hV} \wedge Y)[K_i] 
\rightarrow 
{\varprojlim_i}^s \varinjlim_{H \le V \le_o \GG} 
\pi_*(E^{hV} \wedge Y)[K_i]
\end{equation}
are isomorphisms.  Since the inverse systems are
Mittag-Leffler, we only need to investigate $s = 0$.  By
Lemma~\ref{lem:nilpotent}, the homotopy groups
$\pi_*(E^{hV} \wedge Y)$ and $\pi_*(E^{hH} \wedge Y)$ are finite in each degree.
The map \ref{eq:limcolim} for $s = 0$ is therefore
the composite of the following sequence of
isomorphisms.
\begin{align*}
\varinjlim_{H \le V \le_o \GG} \varprojlim_i \pi_*(E^{hV} \wedge Y)[K_i]
& = \varinjlim_{H \le V \le_o \GG} \pi_*(E^{hV} \wedge Y)[[K]] \\
& \cong \varinjlim_{H \le V \le_o \GG} 
(\pi_*(E^{hV} \wedge Y) \otimes \ZZ[[K]]) \\
& \cong (\varinjlim_{H \le V \le_o \GG} \pi_*(E^{hV} \wedge Y)) 
\otimes \ZZ[[K]] \\
& \cong \pi_*(E^{hH} \wedge Y) \otimes \ZZ[[K]] \\
& \cong \pi_*(E^{hH} \wedge Y)[[K]] \\
& \cong \varprojlim_{i} \pi_*(E^{hH} \wedge Y)[K_i] \\
& \cong \varprojlim_i \varinjlim_{H \le V \le_o \GG} 
\pi_*(E^{hV} \wedge Y)[K_i]
\end{align*}
\qed \end{pf*}

\begin{mycor}\label{cor:Sbaraction}
The action map
$$ \xi : S[[\GG/N]] \wedge E^{hN} \rightarrow E^{hN} $$
of Corollary~\ref{cor:Gaction}
lifts to an action map
$$ \xi : [[\GG/N]]E^{hN} \rightarrow E^{hN}. $$
\end{mycor}

We shall be employing Corollary~\ref{cor:Sbaraction} in
the case of $N = \GG^1$.  The group $\GG/\GG^1 = \ZZ_p$ acts on the 
homotopy fixed point spectrum $\br{S} = E^{h\GG^1}$.
We shall need an additional lemma before proceeding to the proof of
Proposition~\ref{prop:Sbar}.

\begin{mylem}\label{lem:Z_pSES}
Let $X$ be a $p$-complete spectrum, and let $\tau$ be the canonical 
topological generator of the profinite group $\ZZ_p$.  Then the sequence
$$ X[[\ZZ_p]] \xrightarrow{[\tau] - 1} X[[\ZZ_p]] \xrightarrow{\epsilon} X $$
is a cofiber sequence.
\end{mylem}

\begin{pf*}{Proof.}
The map $\epsilon$ is the homotopy inverse limit of the fold maps
$$ X \wedge {\ZZ/p^n}_+ \rightarrow X. $$
Since the composite is null, it will suffice to show that the sequence
induces a short exact sequence on homotopy groups.
Since the inverse system 
$$ \{ \pi_*(X \wedge {\ZZ/p^n}_+) \} = 
\{ \pi_*(X)[\ZZ/p^n] \} $$ 
is Mittag-Leffler, applying $\pi_*$
gives the sequence
$$ \pi_*(X)[[\ZZ_p]] \xrightarrow{[\tau]-1} 
\pi_*(X)[[\ZZ_p]] \rightarrow 
\pi_*(X). $$
This sequence is right exact for any $X$, and
since $\pi_*(X)$ is $p$-complete, this sequence is left exact.
\qed \end{pf*}

\begin{pf*}{Proof of Proposition~\ref{prop:Sbar}}
By Corollary~\ref{cor:FEHEH}, we have a canonical equivalence
\begin{equation*}
\alpha: \br{S}[[\ZZ_p]] \xrightarrow{\simeq} F(\br{S},\br{S})
\end{equation*}
using the fact that $\br{S} = E^{h\GG^1}$ and $\GG/\GG^1 \cong \ZZ_p$.
Using the action of the Morava stabilizer group, given 
an element $\lambda \in \ZZ_p$, we get an induced map 
$$ [\lambda]: \br{S} \rightarrow \br{S}. $$
The effect of the post- and pre-composition maps under $\alpha$ are given
respectively by
\begin{align*}
[\lambda]_* = & L_\lambda^\Delta: \br{S}[[\ZZ_p]] \rightarrow \br{S}[[\ZZ_p]]
\\
[\lambda]^* = & R_\lambda: \br{S}[[\ZZ_p]] \rightarrow \br{S}[[\ZZ_p]]
\end{align*}
where $L_\lambda^\Delta$ is the diagonal left action of $\lambda$ and
$R_\lambda$ us the right action on the $\ZZ_p$ factor.

We shall first prove that $\td{\gamma}$ is an equivalence.
Define a map
$$ \td{\gamma}_1: \br{S}[[\ZZ_p]] \rightarrow F(\br{S}, \br{S}) $$
whose adjoint is the pairing given by the composite
$$ \gamma_1: \br{S}[[\ZZ_p]] \wedge \br{S} 
\xrightarrow{\xi_{\Delta}} \br{S}\wedge\br{S}
\xrightarrow{\mu} \br{S}. $$
Here $\xi_{\Delta}$ is given by allowing the $\ZZ_p$ to act diagonally on
the two $\br{S}$'s (using the action maps of Corollaries~\ref{cor:Sbaraction}
and \ref{cor:Gaction}), and $\mu$ is the product.  We emphasize that
$\td{\gamma}_1$ is a different map from $\alpha$.

Consider the following diagram:
\begin{equation}\label{diag:Sbarcofibers}
\xymatrix{
\br{S}[[\ZZ_p]] \ar[d]^{\td{\gamma}_1}_{\simeq} \ar[r]^{R_\tau-1} 
\ar@{}[dr]|{(1)} &
\br{S}[[\ZZ_p]] \ar[d]^{\td{\gamma}_1}_{\simeq} \ar[r]^{\epsilon} 
\ar@{}[dr]|{(2)} &
\br{S} \ar[d]^{\td{\gamma}}
\\
F(\br{S},\br{S}) \ar[r]_{[\tau]_*-1} &
F(\br{S},\br{S}) \ar[r]_{\delta_*} &
F(\br{S},S^1)
}
\end{equation}
The bottom row is a cofiber sequence, since it arises from the cofiber
sequence~\ref{eq:sbarfiber}.
The top row is a cofiber sequences by Lemma~\ref{lem:Z_pSES}.  
It will follow that
$\td{\gamma}$ is an equivalence once we show that $\td{\gamma}_1$ is an 
equivalence
and we show that squares $(1)$ and $(2)$ of Diagram~\ref{diag:Sbarcofibers}
commute.

Define a map $\beta$ to be the composite
$$ \beta: \br{S}[[\ZZ_p]] \xrightarrow{\Delta} 
\br{S}[[\ZZ_p \times \ZZ_p]] 
\xrightarrow{\xi_L} \br{S}[[\ZZ_p]] $$
where $\Delta$ is the diagonal and $\xi_L$ is given by letting the first
$\ZZ_p$ act on $\br{S}$ through the map of Corollary~\ref{cor:Sbaraction}.  
The map $\beta$ is an isomorphism, with inverse given by the composite
$$ \br{S}[[\ZZ_p]] \xrightarrow{\br{\Delta}} \br{S}[[\ZZ_p\times\ZZ_p]] 
\xrightarrow{\xi_L} \br{S}[[\ZZ_p]] $$
where $\br{\Delta}: \ZZ_p \rightarrow \ZZ_p\times\ZZ_p$ 
is the modified diagonal given by
$$ \br{\Delta}(\lambda) = (\lambda^{-1}, \lambda). $$
(Here, as always, we are using multiplicative notation for the additive
group $\ZZ_p$.)
Then the map $\td{\gamma}_1$ is an equivalence because it is the composite of
the two equivalences
$$ \td{\gamma}_1: \br{S}[[\ZZ_p]] \xrightarrow[\beta]{\simeq} \br{S}[[\ZZ_p]] 
\xrightarrow[\alpha]{\simeq} F(\br{S}, \br{S}). $$

Square~$(1)$ of Diagram~\ref{diag:Sbarcofibers} commutes since it is given
by the following composite, where each of the squares clearly commutes.
$$
\xymatrix{
\br{S}[[\ZZ_p]] \ar[r]^{R_\tau-1} \ar[d]_{\beta} &
\br{S}[[\ZZ_p]] \ar[d]^\beta
\\
\br{S}[[\ZZ_p]] \ar[r]^{L^\Delta_\tau-1} \ar[d]_\alpha &
\br{S}[[\ZZ_p]] \ar[d]^\alpha
\\
F(\br{S},\br{S}) \ar[r]_{[\tau]_*-1} &
F(\br{S},\br{S})
}
$$

Square~$(2)$ of Diagram~\ref{diag:Sbarcofibers} commutes since the adjoint
maps form the following commutative diagram.
$$
\xymatrix{
\br{S}[[\ZZ_p]] \wedge \br{S} \ar[rr]^{\epsilon\wedge 1}
\ar[dr]|{\mu[[Id]]} \ar[d]_{\xi_\Delta} \ar@/_2.5pc/[dd]_{\gamma_1} &&
\br{S}\wedge\br{S} \ar[d]^{\mu} \ar@/^2.5pc/[dd]^{\gamma}
\\
\br{S}\wedge\br{S} \ar@{}[r]|{(3)} \ar[d]_{\mu} &
\br{S}[[\ZZ_p]] \ar[r]^\epsilon \ar[dl]|\xi &
\br{S} \ar[d]^\delta  \ar@{}[dll]|{(4)}
\\
\br{S} \ar[rr]_\delta &&
S^1
}
$$
The map $\xi$ is the action map of $\ZZ_p$ on $\br{S}$ given by
Corollary~\ref{cor:Sbaraction}.  The map $\mu[[Id]]$
is the product on the two factors of $\br{S}$.  The only portions of the
diagram which don't obviously commute are $(3)$ and $(4)$.  Portion~$(3)$
commutes since $\ZZ_p$ acts by maps of $E_\infty$ ring spectra.
Since the composite given in the cofiber sequence
$$ \br{S} \xrightarrow{[\tau]-1} \br{S} \xrightarrow{\delta} S^1 $$
is null, $\delta$ coequalizes the $\ZZ_p$-action.  Therefore, portion~$(4)$
commutes.

We are left with showing that the left-hand square of the diagram in the
statement of the proposition commutes.  It suffices to check on adjoints,
and this is accomplished with the following commuting diagram.
$$
\xymatrix{
\br{S}\wedge S \ar[r]^{1\wedge r} \ar[dd]_{\delta} \ar@{=}[dr] &
\br{S}\wedge \br{S} \ar[d]^\mu 
\\
& \br{S} \ar[d]^\delta 
\\
S^1 \ar@{=}[r] &
S^1
}
$$
\qed \end{pf*}

\subsection{Identification of $I_{K(2)}TMF$}\label{sec:I_KTMF}

In this section we will prove the following proposition.

\begin{myprop}\label{prop:I_KTMF}
The spectra $TMF$ and $TMF_0(2)$ are Gross-Hopkins self dual up to the
suspensions given below.
\begin{align*}
I_KTMF & \simeq \Sigma^{22} TMF \\
I_KTMF_0(2) & \simeq \Sigma^{22} TMF_0(2) \\
I_KE(2) & \simeq \Sigma^{22} E(2) 
\end{align*}
\end{myprop}

For the remainder of this section, we shall work in the $3$-local stable
homotopy category, not localized at $K(2)$. 

Recall from \cite{MahowaldRezk} 
that an $fp$-spectrum is a spectrum whose $\FF_p$-cohomology is finitely
presented over the Steenrod algebra $A$.  
Let $C_n^fX$ be the fiber of the finite 
localization map $X \rightarrow L_n^fX$.
Mahowald and Rezk \cite{MahowaldRezkDuality} studied a dualization functor
$W_n$ defined by
$$ W_n X = IC_n^fX. $$
Their theory indicated that this dual (which we shall refer to as the
Mahowald-Rezk dual) is quite computable for certain $fp$-spectra of type
$n$.  In
particular, they show that both $tmf$ and $BP\bra{2}$ are self-dual.

\begin{myprop}[Mahowald-Rezk \cite{MahowaldRezkDuality}]\label{prop:W_2tmf}
There are equivalences
$$ W_2 tmf \simeq \Sigma^{23} tmf $$
$$ W_2 BP\bra{2} \simeq \Sigma^{23} BP\bra{2}. $$
\end{myprop}

Since $tmf_0(2)$ is equivalent to $BP\bra{2} \vee \Sigma^8 BP\bra{2}$, 
and is
periodic of order $8$, we have the following corollary.

\begin{mycor}
The Mahowald-Rezk dual of $tmf_0(2)$ is given by an equivalence 
$W_2tmf_0(2) \simeq
\Sigma^{23} tmf_0(2)$.
\end{mycor}

We would like to use these computations to identify the Gross-Hopkins duals
of $TMF$ and $TMF_0(2)$.  To this end we investigate the relationship between
the Gross-Hopkins dual and the Mahowald-Rezk dual.
A spectrum $X$ is said to \emph{satisfy the $E(n)$-telescope conjecture} 
if the natural map 
$$ t_n : L_n^f X \rightarrow L_n X $$
is an equivalence.
Our interest in such spectra stems from the following proposition.

\begin{myprop}\label{prop:MRGH}
Suppose that $X$ is a spectrum which satisfies the $E(n)$ and 
$E(n-1)$-telescope
conjectures.  Then there is a cofiber sequence
$$ W_{n-1} X \rightarrow W_n X \rightarrow \Sigma I_n X. $$
\end{myprop}

\begin{pf*}{Proof.}
Applying Verdier's axiom to the composite $X \rightarrow L_n X \rightarrow
L_{n-1}X $ gives a cofiber sequence
$$ \Sigma^{-1} M_nX \rightarrow C_nX \rightarrow C_{n-1}X. $$
Since $X$ satisfies the $E(n)$-telescope conjecture, we may regard this
cofiber sequence as
$$ \Sigma^{-1} M_nX \rightarrow C_n^fX \rightarrow C_{n-1}^fX. $$
The proposition follows after taking the Brown-Comenetz dual.
\qed \end{pf*}

In \cite{MahowaldRezkDuality}, following a suggestion
of Hal Sadofsky, it is pointed out that the
chromatic tower computations of \cite{Ravenel84} imply the following 
proposition.

\begin{myprop}
The spectra $BP\bra{n}$ satisfy the $E(m)$-telescope conjecture for each 
$m$.
\end{myprop}

\begin{pf*}{Proof.}
The spectrum $BP$ satisfies the $E(m)$-telescope conjecture (see 
\cite[7.2]{MahowaldRezk} and \cite[Theorem~6.2]{Ravenel}).
Since the localization functors $L_m^f$ and $L_m$ are smashing, the
$E(m)$-telescope conjecture holds for any spectrum which may be obtained
from $BP$ by means of iterated cofibers and filtered homotopy colimits.  In
particular, the $E(m)$-telescope conjecture holds for the spectrum $BP\bra{n}$.
\qed \end{pf*}

Since $tmf_0(2)$ is equivalent to the wedge $BP\bra{2} \vee \Sigma^8
BP\bra{2}$, one has the following corollary.

\begin{mycor}
The spectrum $tmf_0(2)$ satisfies the $E(m)$-telescope conjecture at the
prime $3$ for every $m$.
\end{mycor}

The following lemma will allow us to 
bootstrap ourselves up from $tmf_0(2)$ to $tmf$.

\begin{lem}
Let $T$ be the $3$-cell complex 
$S^0 \cup_{\alpha_1} e^4 \cup_{\alpha_1} e^8$.
Then there is an equivalence
$$ tmf \wedge T \simeq tmf_0(2). $$
\end{lem}

\begin{pf*}{Proof.}
The spectrum $tmf_0(2)$ is a $tmf$-algebra through the map that forgets the
$\Gamma_0(2)$ structure
$$ \phi_f^*: tmf \rightarrow tmf_0(2). $$
Since the ($3$-primary) homotopy groups of $tmf_0(2) \simeq BP\bra{2} \vee
\Sigma^8 BP\bra{2}$ are concentrated in even degrees, 
the Hurewitz images of the
attaching maps of $T$ in $\pi_*(tmf_0(2))$ are null, and hence the map
$\phi_f^*$ factors to give a map
$$ \td{\phi}_f^* : tmf \wedge T \rightarrow tmf_0(2). $$
The Adams-Novikov $E_2$-term of $tmf \wedge T$, 
and the induced map $\td{\phi}_f^*$ on Adams-Novikov
$E_2$-terms, are easily computed with an
Atiyah-Hirzebruch spectral sequence, from which we deduce that the map
$\td{\phi}_f^*$ induces an isomorphism on the level of Adams-Novikov
$E_2$-terms.  We conclude that the map $\td{\phi}_f^*$ is an equivalence.
\qed \end{pf*}

From this, we deduce the following.

\begin{myprop}
The spectrum $tmf$ satisfies the $E(m)$-telescope conjecture at the prime
$3$ for every $m$.
\end{myprop}

\begin{pf*}{Proof.}
In \cite{Ravenel}, the complex $T$ is used extensively in the following
manner.  One has the cofiber sequences given below.
$$ S^0 \rightarrow T \rightarrow \Sigma^4 C\alpha_1 $$
$$ C\alpha_1 \rightarrow T \rightarrow S^8 $$
Here $C\alpha_1$ is the cofiber of $\alpha_1$.
These splice together to give a $BP\bra{2}$-Adams resolution for $tmf$.
$$
\xymatrix{
tmf \wedge S^0 \ar[d] &
tmf \wedge \Sigma^3 C\alpha_1 \ar[l] \ar[d] &
tmf \wedge S^{10} \ar[l] \ar[d] &
tmf \wedge \Sigma^{13} C\alpha_1 \ar[l] \ar[d] &
\cdots \ar[l]
\\
tmf_0(2) &
\Sigma^3 tmf_0(2) &
\Sigma^{10} tmf_0(2) &
\Sigma^{13} tmf_0(2) 
} $$
Let $X_i$ be the $i^{th}$ term of the tower.  The tower converges in the
sense that there is an equivalence 
$\varprojlim X_i \simeq \ast$.  Let $Y_i$ be given by the
following cofiber sequences
$$ X_i \rightarrow X_0 \rightarrow Y_i. $$
Then by Verdier's axiom there exist cofiber sequences
$$ Y_i \rightarrow Y_{i-1} \rightarrow \Sigma^N tmf_0(2). $$
Since localization respects cofiber sequences, the spectra $Y_i$
inductively satisfy the $E(m)$-telescope conjecture.  We may recover $tmf$
by
$$ tmf \simeq \varprojlim Y_i. $$
Localization does not respect homotopy inverse limits in general,
but this inverse limit is nice, because the composite
$$ S^{10} \rightarrow \Sigma^3 C\alpha_1 \rightarrow S^0 $$
is the element $\beta_1$ of $\pi_{10}(S)$.  Since $\beta_1^6$ is null in
the stable stems, our
tower $X_i$ has the property that $X_{i+14} \rightarrow X_i$ is null.
This makes the tower $Y_i$ pro-equivalent to a constant tower, and thus
localizations commute with this homotopy inverse limit.
We deduce that
$tmf$ satisfies the $E(m)$-telescope conjecture.
\qed \end{pf*}

We recall the following theorem from \cite[Theorem~6.19]{HoveyStrickland}, 
which also appeared without proof in \cite{HopkinsGross}.

\begin{mythm}[\cite{HoveyStrickland}, 
\cite{HopkinsGross}]\label{thm:HopkinsGross}
For any spectrum
$X$, the monochromatic fiber $M_nX$ satisfies
$$ M_nX \simeq M_nL_{K(n)} X. $$
\end{mythm}

\begin{pf*}{Proof of Proposition~\ref{prop:I_KTMF}.}
Let $K$ denote $K(2)$.
We shall prove that $I_KTMF \simeq \Sigma^{-1} L_K W_2tmf$.  
The result
will then follow from Proposition~\ref{prop:W_2tmf}.
The arguments for $tmf_0(2)$ and $BP\bra{2}$ are similar.
Applying $L_{K}$ to the cofiber sequence given by
Proposition~\ref{prop:MRGH}, we get a cofiber sequence
$$ L_K W_{1} tmf \rightarrow L_K W_2 tmf \rightarrow \Sigma I_K tmf. $$
We have an equivalence $L_K tmf \simeq L_K TMF$.  
Theorem~\ref{thm:HopkinsGross} implies that there is an equivalence 
$I_Ktmf \simeq I_KTMF$.  We are left with proving that the map
$$ L_K W_2 tmf \rightarrow \Sigma I_K TMF $$
is an equivalence.  By Lemma~\ref{lem:Bousfield}, it suffices to show that
this map is an isomorphism on $V(1)$-homology.  This is readily seen from
computation.  By Proposition~\ref{prop:W_2tmf} we have
$$ V(1)_*(L_K W_2 tmf) \cong \Sigma^{23} V(1)_*(TMF). $$
Since $L_1V(1) \simeq \ast$, we also have $I_KV(1) \simeq I L_KV(1)$, and
therefore
$$ \Sigma V(1)_*(I_K TMF) \cong \Sigma^7 \pi_*(I(V(1)\wedge TMF)) \cong
\Sigma^{23} V(1)_*TMF $$
The first isomorphism follows from the fact that the CW-complex $V(1)$ is
a Spanier-Whitehead self-dual complex of dimension $6$.  The second
isomorphism follows from the computation 
$$V(1)_*(TMF) \cong \Sigma^{56} 
\Hom(V(1)_*(TMF),\QQ/\ZZ) $$
coupled with the fact that $V(1)_*(TMF)$ is periodic of order $72$.  The
map $L_KW_2tmf \rightarrow \Sigma I_K TMF$ is readily seen to induce an
isomorphism between these $V(1)$-homology groups.
\qed \end{pf*}

\subsection{Identification of $I_{K(2)}^{-1} \wedge TMF$}\label{sec:I_K^TMF}

We now resume the convention that everything is implicitly localized at
$K = K(2)$.
The aim of this section is to prove the following proposition.

\begin{myprop}\label{prop:I_K^TMF}
There are equivalences 
\begin{align*}
I_K \wedge TMF & \simeq \Sigma^{-22} TMF \\
I_K \wedge TMF_0(2) & \simeq \Sigma^{-22} TMF_0(2) \\
I_K \wedge E(2) & \simeq \Sigma^{-22} E(2)
\end{align*}
\end{myprop}

\begin{pf*}{Proof.}
The Gross-Hopkins duality theorem \cite{HopkinsGross}, \cite{Strickland} 
states that there is an isomorphism
$$ E_*(I_K) \cong \Sigma^2 E_*[det] $$
of Morava modules.  Here $[det]$ indicates a twisting with the determinant
character.
There is a homotopy fixed point spectral sequence
\begin{equation}\label{eq:I_K^TMFSS}
H^s(G_{24}; E_t(I_K)) \Rightarrow \pi_{t-s}(TMF \wedge I_K).
\end{equation}
The $E_2$ term of this spectral sequence is thus identified by the
Gross-Hopkins duality theorem.  Since $G_{24}$ is contained in the kernel
of the determinant, we have an isomorphism
$$ H^*(G_{24}; E_*(I_K)) \simeq H^*(G_{24}; \Sigma^2 E_*) $$
The spectral sequence~\ref{eq:I_K^TMFSS} is a spectral sequence of
modules over the spectral sequence of algebras
\begin{equation}\label{eq:TMFSS}
H^s(G_{24}; E_t) \Rightarrow \pi_{t-s}(TMF)
\end{equation}
which computes the homotopy groups of $TMF$.  
The structure of the latter spectral sequence forces the
spectral sequence~\ref{eq:I_K^TMFSS} to be isomorphic to the spectral
sequence~\ref{eq:TMFSS} up to a suspension congruent to $2$ modulo $24$.
Thus there exists a map
$$ S^N \rightarrow TMF \wedge I_K $$
which is detected on the $0$-line of spectral sequence~\ref{eq:I_K^TMFSS}
which extends (using the $TMF$-module structure of $TMF \wedge I_K$) to an
equivalence
$$ \Sigma^N TMF \xrightarrow{\simeq} TMF \wedge I_K. $$
We just need to determine the value of $N$, which is congruent to $2$
modulo $24$ and which is only unique modulo
$72$, the order of periodicity of $TMF$.  
Unfortunately, we are unable to avoid appealing  
to \cite{GoerssHennMahowald} to determine $N$.
This computation was summarized in Section~\ref{sec:piQ(2)V(1)}.  It
implies that
$$ \pi_*(V(1)) \cong \pi_*(\Sigma^{28} IV(1)) \cong \pi_*(\Sigma^{22} I_K
\wedge V(1)) $$
which forces the value of $N$ to be $-22$.

The identification of $I_K \wedge TMF_0(2)$ is similar.  The spectrum
$TMF_0(2)$ has the homotopy fixed point model $E^{hD_8}$.  The group $D_8$ is
in the kernel of the determinant, and the ANSS for $I_K \wedge E^{hD_8}$ is
concentrated on the zero line, so there are no differentials.  Thus one has
an isomorphism
$$ H^*(D_8;E_*I_K) \cong H^*(D_8;\Sigma^2 E_*) $$
which gives an equivalence 
$$ I_K \wedge TMF_0(2) \simeq \Sigma^{2} TMF_0(2) \simeq \Sigma^{-22} TMF_0(2)
$$
using the eightfold periodicity of $TMF_0(2)$.

The case of $E(2)$ is slightly different.  The spectrum $E(2)$ 
is given by the homotopy
fixed point spectrum $E^{hSD_{16}}$, but the subgroup $SD_{16}$ is not
contained in the kernel of the determinant.  In fact, we have an
isomorphism of the restricted determinant 
$$ det \downarrow^\GG_{SD_{16}} \cong \chi $$
where $\chi$ is the non-trivial character of $SD_{16}/D_8$, regarded as an
$SD_{16}$ representation.  Using the fact that there is an isomorphism of
$SD_{16}$-modules
$$ E_* \otimes \ZZ_3^\chi \cong \Sigma^8 E_* $$ 
it follows that there are isomorphisms
$$ H^*(SD_{16};E_*I_K) \cong H^*(SD_{16};\Sigma^2E_*[det]) \cong
H^*(SD_{16}; \Sigma^{10}E_*) $$
which realize to equivalences
$$ I_K \wedge E(2) \simeq \Sigma^{10} E(2) \simeq \Sigma^{-22} E(2) $$
using the $16$-fold periodicity of $E(2)$.
\qed \end{pf*}

\subsection{Proof of Theorem~\ref{thm:DQ}}\label{sec:DQ}

We begin with the an identification of $DTMF$, $DTMF_0(2)$, and $DE(2)$.

\begin{myprop}\label{prop:DTMF}
There are equivalences
\begin{align*}
DTMF & \simeq \Sigma^{44}TMF \\
DTMF_0(2) & \simeq \Sigma^{44}TMF_0(2) \\
DE(2) & \simeq \Sigma^{44} E(2).
\end{align*}
\end{myprop}

\begin{pf*}{Proof.}
Since $I_K$ is invertible, we have for any $X$ \cite{Strickland}
$$ DX \wedge I_K \simeq I_KX. $$
Smashing this equation with $I_K^{-1}$ gives
$$ DX \simeq I_KX \wedge I_K^{-1}. $$
In particular, we may let $X$ be $TMF$, $TMF_0(2)$, or $E(2)$, and then apply 
Propositions~\ref{prop:I_KTMF}
and \ref{prop:I_K^TMF}.  
\qed \end{pf*}

We shall make use of the following lemma, which follows from the
equivalence of fibers and cofibers (up to a suspension) 
in the stable homotopy category.

\begin{mylem}\label{lem:shiftdual}
Suppose $X$ is a spectrum which refines the tower
$$ A_0 \rightarrow A_1 \rightarrow \cdots \rightarrow A_d. $$
Then $DX$ is a spectrum which refines the dual tower
$$ \Sigma^{d} DA_d \rightarrow \Sigma^d DA_{d-1} \rightarrow \cdots
\Sigma^d DA_0. $$
\end{mylem}
Since $Q$ refines a tower of the form
$$ TMF \rightarrow TMF \vee TMF_0(2) \rightarrow TMF_0(2), $$
Theorem~\ref{thm:DQ} follows immediately from
Proposition~\ref{prop:DTMF} and Lemma~\ref{lem:shiftdual}.

\subsection{A spectrum which is half of $\br{S}$}\label{sec:Qbar}

In the next few sections we work our way toward proving
Theorem~\ref{thm:cofiber}.  Proposition~\ref{prop:Sbar} in particular
implies that $\br{S} = E^{h\GG^1}$ is self dual.
In this section we introduce a spectrum $\br{Q}$ which will turn out to be
a Lagrangian for $\br{S}$.

Following \cite[Sec.~4]{GHMR}, let $\chi$ be the the sign representation for
$SD_{16}/D_8$ regarded as an $SD_{16}$-representation. (See
Section~\ref{sec:Morava} for descriptions of these subgroups.)
There is a splitting
$$ TMF_0(2) \simeq E(2) \vee \Sigma^8 E(2)$$
which realizes the splitting
$$ \ZZ_3[[\GG/D_8]] \cong \ZZ_3[[\GG/SD_{16}]] \oplus \ZZ_3^\chi
\uparrow_{SD_{16}}^{\GG} $$
where the $SD_{16}$-representation $\chi$ has been induced up to a
$\GG$-representation.  The first summand is generated by $[\omega]+1 \in
\ZZ_3[[\GG/D_8]]$ and
the second by $[\omega]-1 \in \ZZ_3[[\GG/D_8]]$.  Here, and elsewhere, we
use brackets to denote the group-like elements of our group rings when
confusion may arise as to where the sums are taking place.
Let
$$ \nu : TMF_0(2) \rightarrow \Sigma^8 E(2) $$
be the projection map.

One deficiency of our map $\phi_q^*$ is that the corresponding element
$1+t \in \GG$ (see Section~\ref{sec:Morava}) has norm $det(1+t) =
2$.  Since $\br{S} = E^{h\GG^1}$, where $\GG^1$ is the kernel of the
reduced norm, it would be more natural to consider maps
corresponding to elements of reduced norm $1$.

Let $\sqrt{2} \in \WW$ be the choice of a square root of $2$ which reduces
to $\omega^2$ in $\FF_9$.  Such an
element exists by Hensel's lemma.
We may regard
$\sqrt{2}$ as an element of $\GG$.  
Define a map
$$ \br{\phi}_q^* : TMF \rightarrow TMF_0(2) $$
which in the language of homotopy fixed point spectra, is given by the
composite
$$ E^{hG_{24}} \xrightarrow{\iota_{D_8}^*} E^{hD_8}
\xrightarrow{[(1+t)/\sqrt{2}]} E^{hD_8}. $$
(The Frobenius takes $\sqrt{2}$ to $-\sqrt{2}$, so one may verify that 
$\sqrt{2}$ is an element of the normalizer
$N_\GG D_8$.)

\begin{mydefn}\label{defn:Qbar}
Define $\br{Q}$ to be the homotopy equalizer
$$
\xymatrix{
\br{Q} \ar[r] &
TMF \ar@<1ex>[r]^{\nu \circ \phi_f^*} \ar@<-1ex>[r]_{\nu \circ \br{\phi}_q^*} &
\Sigma^8 E(2).
} $$
\end{mydefn}

The following lemma implies that we may transfer our understanding of the
effect of 
$\phi_q^*$ on $V(1)$-homology to that of $\br{\phi}^*_q$.

\begin{mylem}\label{lem:phisub}
The maps
$$ V(1)_*\phi^*_q, V(1)_*\br{\phi}^*_q : V(1)_*(TMF) \rightarrow
V(1)_*(TMF_0(2)) $$
are identical.
\end{mylem}

\begin{pf*}{Proof.}
The map $\br{\phi}^*_q$ differs from $\phi_q^*$ by a factor of $\sqrt{2}
\in \WW$.  Since $\sqrt{2}$ reduces to $\omega^2 \in \FF_9$, and
$V(1)_*TMF_0(2)$ is concentrated in degrees congruent to $0$ modulo $8$, the
effect of $\sqrt{2}$ on $V(1)_*TMF_0(2)$ is multiplication by $\omega^8 = 1$.
\qed \end{pf*}

Let $e_8$ be the generator in degree $8$ of $\pi_*(\Sigma^8 E(2))$.  Then
our computations in Section~\ref{sec:mapscomp} combined with
Lemma~\ref{lem:phisub} imply that when we compute
the effects of the map on Adams-Novikov $E_2$-terms
$$ \nu_* \circ (\br{\phi}_q^*-\phi_f^*): E_2(V(1) \wedge TMF) \rightarrow
E_2(V(1) \wedge \Sigma^8 E(2)) $$
we have
$$ \nu_* \circ (\br{\phi}_q^*-\phi_f^*)(q_4^k) = 
\begin{cases}
-v_2^{(k-1)/2} e_8, & \text{if $k$ is odd} \\
0 & \text{if $k$ is even}
\end{cases}
$$
From these formulas one easily computes the ANSS $E_2$-term $E_2(\br{Q})$.
The ANSS differentials are induced from the differentials in the ANSS for
$TMF$ just as in Proposition~\ref{prop:Q(2)V(1)ANSSdiffs}.
One finds that (using the patterns of homotopy groups described in
Section~\ref{sec:Q(2)V(1)}) we have an isomorphism
\begin{equation}\label{eq:V(1)Qbar}
V(1)_*(\br{Q}) = (B \oplus C)\otimes P[v_2^{\pm9}].
\end{equation}
This should be compared to the computation of the 
$V(1)$-homology of $\br{S}$ given in \cite{GoerssHennMahowald}.
\begin{equation}\label{eq:V(1)Sbar} 
V(1)_*(\br{S}) = 
((B \oplus C) \oplus \Sigma^{29} (B
\oplus C)^{\vee})\otimes P[v_2^{\pm 9}]
\end{equation}
We see that
$V(1)_*(\br{Q})$ is an additive summand of both $V(1)_*(\br{S})$ and
$V(1)_*(Q)$, and that $V(1)_*(\br{S})$ is the direct sum of a copy of
$V(1)_*(\br{Q})$ and a shifted copy of the dual of $V(1)_*(\br{Q})$.

\subsection{Recollections from \cite{GHMR}}\label{sec:GHMR}

In \cite[Prop.~2.6]{GHMR} the mapping spectra $F(E^{hF_1},E^{hF_2})$ are 
described for $F_i$
closed subgroups of $\GG$, with $F_2$ finite.  The authors of \cite{GHMR} prove
the following proposition.

\begin{myprop}[Goerss-Henn-Mahowald-Rezk, \cite{GHMR}]\label{prop:GHMRmapping}
We have
$$ F(E^{hF_1},E^{hF_2}) \simeq \prod_{x \in F_2\backslash \GG / F_1} 
E^{hF_x} $$
where $F_x$ is the finite subgroup $F_2 \cap x F_1 x^{-1}$.
\end{myprop}

The product in Proposition~\ref{prop:GHMRmapping} must be properly
interpreted as an appropriate homotopy inverse limit, using the profinite
structure of the double coset space $F_2\backslash \GG / F_1$.

We are primarily concerned with the following spectra.
\begin{align*}
E^{hG_{24}} & = TMF \\
E^{hD_8} & = TMF_0(2) \\
E^{hSD_{16}} & = E(2)
\end{align*}

The homotopy groups of $E^{hF}$ are computed in \cite{GHMR} for $F$ finite.
We shall use the following
observations about $\pi_*E^{hF}$ for $F$ finite.

\begin{mylem}\label{lem:EhF24}
Suppose that one of the $F_i$ is $G_{24}$ and the other is either $G_{24}$,
$D_8$, or $SD_{16}$.  Then the homotopy of each of the summands
$\pi_*(E^{hF_x})$ described in Proposition~\ref{prop:GHMRmapping} is
concentrated in degrees congruent to $0$ modulo $4$ and
$1,3,10,13,27,30$ modulo $36$.
\end{mylem}

\begin{pf*}{Proof.}
The group $F_x = F_2 \cap x F_1 x^{-1}$ is contained in a conjugate of
the subgroup $G_{24}$, and contains the central element $-1$, since each of
the subgroups $G_{24}$, $D_8$, and $SD_{16}$ contains the element $-1$.  We
conclude that either: (1) the group $F_x$ contains a conjugate of the cyclic
group $C_6 = \bra{s, -1}$, or (2) the order of the group $F_x$ is prime to
$3$, and $F_x$ contains the element $-1$.  In
either case $F_x$ contains the element $-1$.
The ring of
invariants $E_*^{F_x}$ is concentrated in degrees congruent to $0$ modulo
$4$, since this is true of the invariants $E_*^{\pm 1}$
\cite[Cor.~3.15]{GHMR}.  If we are in case
(2), then there is no higher cohomology, and the homotopy fixed point
spectral sequence collapses to give
$$ \pi_*(E^{hF_x}) = E_*^{F_x}. $$
If we are in case (1), then the computations of 
\cite[Thm~3.10, Rmk.~3.12]{GHMR}
indicate that elements of $\pi_*(E^{hF_x})$
arising in the homotopy fixed
point spectral sequence from group
cohomology elements of cohomological degree greater than $0$ 
may only lie in degrees 
congruent to $1$, $3$, $10$, $13$, $27$, or $30$ modulo $36$.
\qed \end{pf*}

\begin{mylem}\label{lem:EhF4}
Suppose that one of the $F_i$ is either $D_8$ or $SD_{16}$ and the other is
either $G_{24}$, $D_8$, or $SD_{16}$.
Then the homotopy of each of the summands
$\pi_*(E^{hF_x})$ described in Proposition~\ref{prop:GHMRmapping} is
concentrated in degrees congruent to $0$ modulo $4$.
\end{mylem}

\begin{pf*}{Proof.}
Since one of the subgroups $F_i$ has order prime to $3$, and both of the
subgroups $F_i$ contain the central element $-1$, we are in case (2) of the
proof of Lemma~\ref{lem:EhF24}
\qed \end{pf*}

\subsection{A Lagrangian decomposition of
$\br{S}$}\label{sec:Sbarlagrangian}

Our calculations of the $V(1)$-homology of $\br{S} = E^{h\GG^1}$
(\ref{eq:V(1)Sbar}) and 
$\br{Q}$ (\ref{eq:V(1)Qbar}) suggest the
following Lagrangian decomposition.

\begin{myprop}\label{prop:Sbarlagrangian}
There is a map $\br{\eta}$ such that the following is a fiber sequence
$$ \Sigma D\br{Q} \xrightarrow{D\br{\eta}} \br{S} \xrightarrow{\br{\eta}}
\br{Q}.
$$
The map $D\br{\eta}$ is the dual of the map $\br{\eta}$, using the
equivalence $D\br{S} \simeq \Sigma^{-1} \br{S}$ given by
Proposition~\ref{prop:Sbar}.
\end{myprop}

\begin{myrem}
Proposition~\ref{prop:Sbarlagrangian} provides an alternative construction of
the resolution of $\br{S}$ given in \cite{GHMR}.
\end{myrem}

The inclusion $\iota_{G_{24}}: G_{24} \hookrightarrow \GG^1$ induces a map
$$ \iota_{G_{24}}^* : \br{S} \rightarrow TMF. $$
The spectrum $\br{Q}$ is the homotopy equalizer of the maps
$ \nu \circ \br{\phi}_q^*$ and $\nu \circ \phi_f^*$
(Definition~\ref{defn:Qbar})
where $\br{\phi}_q^*$ and $\phi_f^*$ correspond to the elements
$(t+1)/\sqrt{2}$ and $1$ of $\GG$, respectively.  
Both of these elements have norm $1$.
Therefore the two composites
$$ 
\xymatrix{
\br{S} \ar[r]^{\iota_{G_{24}}^*} & 
TMF \ar@<1ex>[r]^{\br{\phi}_q^*} \ar@<-1ex>[r]_{\phi_f^*} &
TMF_0(2) \ar[r]^\nu &
\Sigma^8 E(2) 
}
$$
actually agree.  Thus there is a lift of $\iota_{G_{24}}^*$ 
to a map
$$ \br{S} \rightarrow \br{Q}. $$
This is the map $\br{\eta}$.

\begin{myrem}
The lift $\br{\eta}$ of $\iota^*_{G_{24}}$ is unique, since
Proposition~\ref{prop:GHMRmapping} implies there are no nontrivial maps
$\br{S} \rightarrow \Sigma^7 E(2)$.
\end{myrem}

\begin{mylem}\label{lem:compositenull}
The composite
$$ \Sigma D\br{Q} \xrightarrow{D\br{\eta}} \br{S} \xrightarrow{\br{\eta}}
\br{Q} $$
is null.
\end{mylem}

\begin{pf*}{Proof.}
We  shall prove that
$$ [\Sigma D \br{Q}, \br{Q}] = 0. $$
Since $\br{Q}$ is built from $TMF$ and $\Sigma^7E(2)$, by
Proposition~\ref{prop:DTMF}, $\Sigma D\br{Q}$ is built from
$\Sigma^{45}TMF$ and $\Sigma^{38} E(2)$.  Lemma~\ref{lem:EhF24} implies that
the following groups are zero.
\begin{alignat*}{2}
\pi_{45}(F(TMF,TMF)) & = 0 & \qquad \pi_{38}(F(TMF, E(2))) & = 0 \\
\pi_{38}(F(E(2),TMF)) & = 0 & \qquad \pi_{31}(F(E(2),E(2))) & = 0
\end{alignat*}
It follows that there are no essential maps $\Sigma D\br{Q} \rightarrow
\br{Q}$.
\qed \end{pf*}

\begin{pf*}{Proof of Proposition~\ref{prop:Sbarlagrangian}}
Let $F$ be the fiber of the map $\br{\eta}$.  
Lemma~\ref{lem:compositenull} implies that there exists a lift $f$ making the
following diagram commute.
$$
\xymatrix{
F \ar[r]^\iota &
\br{S} \ar[r]^{\br{\eta}} &
\br{Q} &
\\
& \Sigma D\br{Q} \ar@{.>}[ul]^f \ar[u]_{D\br{\eta}} 
} $$
Consider the maps $V(1)_*\iota$ and $V(1)_*D\br{\eta}$ on
$V(1)$-homology.  The map $V(1)_*\iota$ is an isomorphism onto $\ker
V(1)_*\br{\eta}$ since $V(1)_*\br{\eta}$ is surjective.  The map 
$V(1)_*D\br{\eta}$ is seen to be an isomorphism onto $\ker V(1)_*\br{\eta}$
from our explicit knowledge of the $V(1)$-homology groups.  The effect of
the map 
$V(1)_*D\br{\eta}$ on $V(1)$-homology is determined from $V(1)_*\br{\eta}$,
since $V(1)_*D\br{\eta}$ is just the Pontryagin dual of $V(1)_*\br{\eta}$.
Since both $V(1)_*\iota$ and $V(1)_*D\br{\eta}$ induce isomorphisms 
onto their
images, $V(1)_*f$ must be an isomorphism.  By Lemma~\ref{lem:Bousfield}
the map $f$ must therefore be an equivalence.
\qed \end{pf*}

\subsection{Building $Q$ from $\br{Q}$}\label{sec:Qdecomp}

In the first part of this section we will produce a map
$$ Q \xrightarrow{r_Q} \br{Q} $$
which will turn out to 
be equivalent to the unit map $Q = Q \wedge S \xrightarrow{1\wedge r} Q \wedge
\br{S}$.  The map $r_Q$ will be produced by the methods of \cite[Sec.~4]{GHMR}.
Specifically, we will prove the following lemma.

\begin{mylem}\label{lem:ZZGGmods}
There exists a continuous map $\rho$ of
$\ZZ_3[[\GG]]$-modules making the following diagram commute
\begin{equation}\label{diag:ZGmods}
\xymatrix@C+2em{ 
\ZZ_3[[\GG/G_{24}]] \ar@{<-}[r]^{[(t+1)/\sqrt{2}]-1} \ar@{=}[d] & 
\ZZ_3[[\GG/D_8]] \ar@{<-}[r]^{\nu} & 
\ZZ_3^\chi \uparrow_{SD_{16}}^{\GG} \ar@{.>}[d]^{\rho}
\\
\ZZ_3[[\GG/G_{24}]] \ar@{<-}[rr]_-{([t+1]-1)\oplus([2]-1)} &&
\ZZ_3[[\GG/D_8]] \oplus \ZZ_3[[\GG/G_{24}]]
}
\end{equation}
\end{mylem}

We postpone the proof of Lemma~\ref{lem:ZZGGmods}.

\begin{mylem}
The map $\rho$ induces a map $\rho^*$ making the following diagram
commute.
$$
\xymatrix{
TMF \ar[r]^{\nu\circ(\br{\phi}^*_q-\phi_f^*)} \ar@{=}[d] &
\Sigma^8 E(2) 
\\
TMF \ar[r]_-{d_0-d_1} & 
TMF_0(2) \vee TMF \ar[u]^{\rho^*}
} $$
\end{mylem}

\begin{pf*}{Proof.}
By Proposition~\ref{prop:GHMRmapping}, for $H$ finite, 
the mapping space $F(E^{hH}, E(2))$
is equivalent to a homotopy inverse limit of products of homotopy fixed
point spectra $E^{hG}$ with $G$ contained in $SD_{16}$.  In particular, the
groups $G$ have order prime to $3$, and thus the ANSS for the mapping space
$F(E^{hH}, E(2))$ is concentrated in the zero line.

Therefore, the edge homomorphism gives a sequence of isomorphisms
\begin{align*}
[E^{hH}, \Sigma^{8}E(2)]
& \cong \Hom_{\ZZ_3[SD_{16}]}^c(\ZZ_3,\pi_8(E)[[\GG/H]]) \\
& \cong 
\Hom_{\ZZ_3[SD_{16}]}^c(\ZZ_3,\ZZ_3^{\chi^{-1}}\widehat{\otimes} \ZZ_3[[\GG/H]] 
) \\
& \cong \Hom_{\ZZ_3[SD_{16}]}^c(\ZZ_3^\chi,\ZZ_3[[\GG/H]]) \\
& \cong 
\Hom_{\ZZ_3[[\GG]]}^c(\ZZ_3^\chi\uparrow_{SD_{16}}^{\GG},\ZZ_3[[\GG/H]]).
\end{align*}
It follows that the map $\rho$ gives rise to the desired map, and the
commutativity of the square in the statement of the lemma follows from the
commutativity of Diagram~\ref{diag:ZGmods}.
\qed \end{pf*}

The map of towers 
$$
\xymatrix{
TMF \ar[r] & 
\Sigma^{8} E(2) \ar[r] &
\ast 
\\
TMF \ar@{=}[u] \ar[r] &
TMF_0(2) \vee TMF \ar[u]^{\rho^*} \ar[r] &
TMF_0(2) \ar[u]
}
$$
induces a map
$$ r_Q: Q \rightarrow \br{Q}. $$ 
The importance of this map for us is encoded in the following lemma and its
corollary.

\begin{mylem}\label{lem:xi}
There is an equivalence $\xi$ making the following diagram commute.
$$
\xymatrix{
Q \ar[r]^-{1\wedge \eta} \ar[dr]_{r_Q} & 
Q \wedge \br{S} \ar@{.>}[d]^{\xi}_{\simeq}
\\
& \br{Q}
} $$
\end{mylem}

The proof of Lemma~\ref{lem:xi} is postponed.

\begin{mycor}\label{cor:Sbarcofiber}
There is a fiber sequence
$$ \Sigma^{-1} \br{Q} \rightarrow Q \xrightarrow{r_Q} \br{Q}. $$
\end{mycor}

\begin{pf*}{Proof.}
Simply smash the fiber sequence
$$ \Sigma^{-1} \br{S} \xrightarrow{\delta} S \xrightarrow{r} \br{S} $$
with $Q$, and apply Lemma~\ref{lem:xi}.
\qed \end{pf*}

In order to prove Lemma~\ref{lem:ZZGGmods}, 
we shall need the following technical result.

\begin{mylem}\label{lem:rightexact}
The sequence
$$ \ZZ_3[[\GG/D_8]] \oplus \ZZ_3[[\GG/G_{24}]] 
\xrightarrow{f} \ZZ_3[[\GG/G_{24}]]
\xrightarrow{\epsilon} \ZZ_3 \rightarrow 0
$$
is exact, where $\epsilon$ is the augmentation and $f$ is the map 
$([t+1]-1)\oplus([2]-1)$.
\end{mylem}

\begin{pf*}{Proof.}
The composite $\epsilon \circ f$ is clearly zero.  
Let $N \subset \ZZ_3[[\GG/G_{24}]]$ be the kernel of $\epsilon$.
We must show that $f$ surjects onto $N$.
Lemma~4.3 of \cite{GHMR}
implies that it suffices to show that the reduced map
$$ \FF_3 \otimes f : \FF_3 \otimes_{\ZZ_3[[\GG]]} 
\left( \ZZ_3[[\GG/D_8]] \oplus \ZZ_3[[\GG/G_{24}]] \right) \rightarrow 
\FF_3 \otimes_{\ZZ_3[[\GG]]} N $$
is surjective.  This is equivalent to showing that
$$ f^* : \ext^0_{\ZZ_3[[\GG]]}(N,\FF_3) \rightarrow 
\ext^0_{\ZZ_3[[\GG]]}(\ZZ_3[[\GG/D_8]] \oplus \ZZ_3[[\GG/G_{24}]],\FF_3) $$
is injective.  The group $\ext^0(N,\FF_3)$, and the map $f^*$, may be easily
deduced from our computations of the ANSS $E_2$-terms for $V(1)$-homology,
from which it is explicitly seen that indeed $f^*$ is injective.
\qed \end{pf*}

\begin{pf*}{Proof of Lemma~\ref{lem:ZZGGmods}.}
For any $\ZZ_3[[G]]$ module $M$, producing a map
$$ \ZZ_3^\chi \uparrow^\GG_{SD_{16}} \rightarrow M $$
is equivalent to specifying an element of $M$ on which $SD_{16}$ acts on
with the sign representation.  Let $x \in \ZZ_3[[\GG/G_{24}]]$ be the
element corresponding to the map
$$ ([(t+1)/\sqrt{2}]-1)\circ \nu : \ZZ_3^\chi 
\uparrow^\GG_{SD_{16}} \rightarrow
\ZZ_3[[\GG/G_{24}]]. $$
The element $x$ is in the kernel of $\epsilon$, hence by
Lemma~\ref{lem:rightexact} it lifts to an element $\td{x}$ of 
$\ZZ_3[[\GG/D_8]] \oplus \ZZ_3[[\GG/G_{24}]]$.  While it is not necessarily
true that $\td{x}$ generates a copy of $\chi$, the weighted average
$$ y = \frac{1}{16}\sum_{g \in SD_{16}} \chi(g) [g]\td{x} $$
does have the property that $\ZZ_3y \cong \ZZ_3^\chi$.  
The element $y$ corresponds
to the map $\rho$.
\qed \end{pf*}

\begin{pf*}{Proof of Lemma~\ref{lem:xi}}
We first must define $\xi$.
Note that since $\nu
\circ \br{\phi}^*_q$ and $\nu \circ \phi_f^*$ are maps of $\br{S}$-modules,
the homotopy equalizer $\br{Q}$ is a $\br{S}$-module.  Therefore there is
a right action map
$$ \mu: \br{Q} \wedge \br{S} \rightarrow \br{Q}. $$
The map $\xi$ is defined to be the composite
$$ \xi: Q \xrightarrow{r_Q \wedge r} \br{Q} \wedge \br{S} \xrightarrow{\mu}
\br{Q} $$
where $r$ is the unit $S \rightarrow \br{S}$.

Recall that, in the notation of Section~\ref{sec:piQ(2)V(1)}, we have 
$$ V(1)_*(Q) = (B \oplus C)\otimes P[v_2^{\pm 9}]\otimes E[\zeta] $$
$$ V(1)_*(\br{Q}) = (B \oplus C)\otimes P[v_2^{\pm 9} ]. $$
The definition of $r_Q$ implies that diagram below commutes.
$$
\xymatrix{
Q \ar[r]^{r_Q} \ar[d] &
\br{Q} \ar[d] 
\\
TMF \ar@{=}[r] &
TMF
} $$
Since these are maps of $S$-modules, the only possibility is for the
induced map $V(1)_*r_Q$ on $V(1)$-homology is to be the quotient by the ideal
generated by $\zeta$.

In Appendix~A, we prove Theorem~\ref{thm:Sbarcell}, which 
says that $\br{S}$ has a cell decomposition 
$$ S^0 \cup_{\zeta} e^0 \cup_{\zeta} e^0 \cdots .$$
It follows that 
$V(1)_*(Q \wedge \br{S})$ is the quotient of $V(1)_*(Q)$ by the ideal
generated by $\zeta$, and 
the map
$$ (1 \wedge r)_* : \pi_*(Q) \rightarrow \pi_*(Q \wedge \br{S})  $$
is the quotient map.

Our map $\xi$ is easily shown to make the diagram of in the statement of
Lemma~\ref{lem:xi} commute, and since $V(1)_*(r_Q)$ and $V(1)_*(1 \wedge r)$
are both surjections onto the same image, $V(1)_*(\xi)$ must be an
isomorphism.  Lemma~\ref{lem:Bousfield} implies that $\xi$ is an
equivalence.
\qed \end{pf*}

\subsection{Proof of Theorem~\ref{thm:cofiber}}\label{sec:cofiberproof}

In this section we piece together our Lagrangians $\br{Q}$ to prove 
Theorem~\ref{thm:cofiber}.  We first observe that the sequence of
Theorem~\ref{thm:cofiber} looks like a fiber sequence on $V(1)$-homology.

\begin{mylem}\label{lem:V(1)SES}
Using the notation of Section~\ref{sec:piQ(2)V(1)}, there is an isomorphism
of short exact sequences.
$$
\xymatrix{
V(1)_*DQ \ar@{=}[d] \ar[r]^{D\eta} &
V(1)_* S \ar@{=}[d] \ar[r]^\eta &
V(1)_* Q \ar@{=}[d]  
\\
\Sigma^{29}(B \oplus C)^{\vee}\otimes R \ar[r] &
\left(  
\genfrac{}{}{0pt}{}{(B \oplus C) \oplus }{\Sigma^{29}(B 
\oplus C)^\vee}
\right) \otimes R \ar[r] &
(B \oplus C) \otimes R
} $$
where $R$ is the ring $P[v_2^{\pm 9}] \otimes E[\zeta]$ and the maps in the
bottom row are the obvious inclusions and projections.
\end{mylem}

\begin{pf*}{Proof.}
The map $\eta$ was essentially computed in Section~\ref{sec:piQ(2)V(1)}.
The map $D\eta$ is, up to suspension, just the Pontryagin dual of the map
$\eta$.
\qed \end{pf*}

We shall need the following lemma.

\begin{mylem}\label{lem:DQ^Sbar}
There is an equivalence $DQ\wedge \br{S} \simeq \Sigma D\br{Q}$.
\end{mylem}

\begin{pf*}{Proof.}
Although $\br{S}$ is not dualizable, we will nevertheless show that the
natural map
$$ \wedge: DQ \wedge D\br{S} \rightarrow D(Q\wedge \br{S})$$
is an equivalence.  The equivalence of the statement of the lemma is then
the composite
$$ DQ \wedge \Sigma^{-1}\br{S} \xrightarrow[\simeq]{1\wedge \td{\gamma}} DQ
\wedge D\br{S} \xrightarrow{\wedge} D(Q\wedge\br{S})
\xrightarrow[\simeq]{D\xi^{-1}} D\br{Q} $$
where $\td{\gamma}$ is the equivalence given in Proposition~\ref{prop:Sbar}, 
and $\xi$ is the equivalence given in Lemma~\ref{lem:xi}.

Let $r: S \rightarrow \br{S}$ be the unit.
Consider the following commutative diagram.
\begin{equation}\label{diag:Dr}
\xymatrix{
DQ \wedge D\br{S} \ar[r]^{1\wedge Dr} \ar[d]_{\wedge} &
DQ \wedge DS \ar[d]^{\wedge}_{\simeq}
\\
D(Q \wedge \br{S}) \ar[r]_{D(1\wedge r)} &
DQ
}
\end{equation}
We claim that when we apply $V(1)$-homology to Diagram~\ref{diag:Dr}, 
we get the following, using the notation
of Section~\ref{sec:piQ(2)V(1)}
\begin{equation}\label{diag:V(1)effect}
\xymatrix@C+2em{
\Sigma^{28}(B \oplus C)^\vee \otimes P \quad 
\ar@{>->}[r]^{V(1)_*1\wedge Dr}_{\cdot \zeta} \ar[d]_{V(1)_*\wedge} &
\Sigma^{29}(B \oplus C)^\vee \otimes R
\ar@{=}[d]^{V(1)_*\wedge}
\\
\Sigma^{28}(B \oplus C)^\vee \otimes P \quad
\ar@{>->}[r]_{V(1)_*D(1\wedge r)}^{\cdot \zeta} &
\Sigma^{29}(B \oplus C)^\vee \otimes R
}
\end{equation}
where $P = P[v_2^{\pm9}]$ and $R = P \otimes E[\zeta]$.
We are claiming that the top and bottom maps are the inclusions
given by multiplication
by $\zeta$.  We see that the left hand map $V(1)_*\wedge$ must therefore be
an isomorphism.  Therefore, by Lemma~\ref{lem:Bousfield}, it is an
equivalence, and the lemma is proven.

We are left with showing that the top and bottom maps in
Diagram~\ref{diag:V(1)effect}
behave as claimed on $V(1)$-homology.  The cellular model of $\br{S}$ given
in Theorem~\ref{thm:Sbarcell} implies that there is an the following 
isomorphism of short
exact sequences.
\begin{equation}\label{diag:SESs}
\xymatrix@C-.5em{
0 \ar[r] &
V(1)_*(DQ\wedge \Sigma^{-1}\br{S}) \ar[r]^{1\wedge \delta} \ar@{=}[d] &
V(1)_*(DQ\wedge S) \ar[r]^{1\wedge r} \ar@{=}[d] &
V(1)_*(DQ\wedge \br{S}) \ar[r] \ar@{=}[d] &
0
\\
0 \ar[r] &
\Sigma^{28}(B\oplus C)^\vee \otimes P \ar[r]_{\cdot \zeta} &
\Sigma^{29}(B \oplus C)^\vee \otimes R \ar[r] &
\Sigma^{29}(B \oplus C)^\vee \otimes P \ar[r] &
0
}
\end{equation}

Smashing the left hand square of the diagram of Proposition~\ref{prop:Sbar}
with $DQ$, we get a commutative diagram
$$
\xymatrix{
DQ \wedge \Sigma^{-1}\br{S} \ar[r]^{1\wedge \delta}
\ar[d]_{1\wedge\td{\gamma}}^{\simeq} &
DQ\wedge S \ar@{=}[d]
\\
DQ\wedge D\br{S} \ar[r]_{1\wedge Dr} &
DQ \wedge S
}
$$
from which it follows that $1 \wedge Dr$ behaves as advertised on
$V(1)$-homology.

Taking the Spanier-Whitehead dual of the diagram of Lemma~\ref{lem:xi}
gives the following commutative diagram.
$$
\xymatrix{
D(Q \wedge \br{S}) \ar[r]^-{D(1 \wedge r)} \ar[d]_{D\xi}^{\simeq} &
DQ \ar@{=}[d] 
\\
D\br{Q} \ar[r]_{Dr_Q} &
DQ
} $$
Gross-Hopkins duality implies that $V(1)$-cohomology is, 
up to a shift in degree, the Pontryagin dual of $V(1)$-homology.
Therefore, the description of $V(1)_*(r_Q)$ given in
the proof of Lemma~\ref{lem:xi} dualizes to give the desired behavior of 
$V(1)_*D(1\wedge r)$ in Diagram~\ref{diag:V(1)effect}.
\qed \end{pf*}

We are now ready to prove a weaker version of Theorem~\ref{thm:cofiber}
where we have smashed everything with $\br{S}$.

\begin{mylem}\label{lem:Sbarsmashcofiber}
The sequence
$$ DQ \wedge \br{S} \xrightarrow{D\eta \wedge 1} S \wedge \br{S} 
\xrightarrow{\eta \wedge 1} 
Q \wedge \br{S} $$
is a fiber sequence.
\end{mylem}

\begin{pf*}{Proof.}
The cellular description of $\br{S}$ given in Theorem~\ref{thm:Sbarcell}
implies that
on $V(1)$-homology, smashing with $\br{S}$ corresponds to modding out by
the ideal generated by $\zeta$.  Using Lemma~\ref{lem:V(1)SES}, we see that
we have an isomorphism of short exact sequences
$$
\xymatrix{
V(1)_*DQ \wedge \br{S} \ar@{=}[d] \ar[r]^{D\eta} &
V(1)_* \br{S} \ar@{=}[d] \ar[r]^\eta &
V(1)_* Q\wedge \br{S} \ar@{=}[d]  
\\
\Sigma^{29}(B \oplus C)^{\vee}\otimes P \ar[r] &
\left(  
\genfrac{}{}{0pt}{}{(B \oplus C) \oplus }{\Sigma^{29}(B 
\oplus C)^\vee}
\right) \otimes P \ar[r] &
(B \oplus C) \otimes P
} $$
where $P = P[v_2^{\pm 9}]$.

By Lemma~\ref{lem:xi}, we have an equivalence $Q \wedge \br{S} \simeq
\br{Q}$, and by Lemma~\ref{lem:DQ^Sbar}
we have an equivalence
$DQ \wedge \br{S} \simeq \Sigma D\br{Q}$.
In the proof of Lemma~\ref{lem:compositenull}, it was proven that 
$$ [\Sigma D \br{Q}, \br{Q} ] = 0. $$
Thus the composite
$$ DQ \wedge \br{S} \xrightarrow{D\eta \wedge 1} \br{S} \xrightarrow{\eta
\wedge 1} 
Q \wedge \br{S} $$
is null.  
The induced map from $DQ \wedge \br{S}$ to the fiber of $\eta
\wedge 1$ is seen to be a $V(1)$-homology equivalence, hence by
Lemma~\ref{lem:Bousfield}, it is an equivalence.
\qed \end{pf*} 

We will now complete the proof of Theorem~\ref{thm:cofiber}.  The reason that
Theorem~\ref{thm:cofiber} is more difficult to prove than
Proposition~\ref{prop:Sbarlagrangian} is that the composite
$$ DQ \xrightarrow{D\eta} S \xrightarrow{\eta} Q$$
cannot be shown to be null by dimensional considerations alone.  However, we
have shown in Lemma~\ref{lem:Sbarsmashcofiber} that after smashing the
above sequence with $\br{S}$ we get a cofiber sequence.

\begin{mylem}\label{lem:fibermap}
The map $\tau - 1: \br{S} \rightarrow
\br{S}$, induces a map of fiber sequences
$$
\xymatrix@C+2em{
\Sigma^{-1}Q \wedge \br{S} \ar[r]^{\delta'} \ar[d]|{1 \wedge (\tau-1)} &
DQ \wedge \br{S} \ar[d]|{1 \wedge (\tau-1)} \ar[r]^{D\eta\wedge 1} &
S \wedge \br{S} \ar[d]|{1 \wedge(\tau-1)} \ar[r]^{\eta \wedge 1} &
Q \wedge \br{S} \ar[d]|{1 \wedge (\tau-1)} 
\\
\Sigma^{-1}Q \wedge \br{S} \ar[r]^{\delta'} &
DQ \wedge \br{S} \ar[r]^{D\eta\wedge 1} &
S \wedge \br{S} \ar[r]^{\eta \wedge 1} &
Q \wedge \br{S}
} $$
where $\delta'$ is the completion of the sequence of
Lemma~\ref{lem:Sbarsmashcofiber} to a fiber seqeuence.
\end{mylem}

\begin{pf*}{Proof.}
We need to show that the map $\delta'$ commutes with the map $1 \wedge
(\tau-1)$.  Lemma~\ref{lem:xi}, Lemma~\ref{lem:DQ^Sbar}, and
Proposition~\ref{prop:DTMF} combine to show that there are fiber sequences
\begin{gather*}
\Sigma^6 E(2) \rightarrow 
\Sigma^{-1} Q \wedge \br{S} \rightarrow \Sigma^{-1} TMF \rightarrow 
\Sigma^7 E(2)
\\
\Sigma^{5} E(2) \rightarrow \Sigma^{45} TMF \rightarrow DQ \wedge \br{S}
\rightarrow \Sigma^{6} E(2).
\end{gather*}
By Lemma~\ref{lem:EhF24}, the groups 
$$ [\Sigma^{-1}TMF, \Sigma^{45}TMF], \quad
[\Sigma^{-1}TMF, \Sigma^{6}E(2)], \quad \mathrm{and} \quad
[\Sigma^{6}E(2), \Sigma^{45}TMF]
$$
are all zero, so we conclude that the induced map
$$ [\Sigma^{-1}Q \wedge \br{S}, DQ \wedge \br{S}] 
\rightarrow [\Sigma^{6}E(2), \Sigma^{6}E(2)]
$$
is a monomorphism.  The topological generator 
$$ \tau \in \ZZ_p = \ZZ_p^\times/\FF_p^\times = \GG/\GG^1 $$
may be lifted to a perfect square $a^2$ in $\ZZ_p^\times$, hence to a
central element $a$ in $\GG$.  The lemma therefore follows from the fact
that the map $[a]-1$ is central in endomorphism ring
$$ [E(2),E(2)] = (E(2)_*[[\GG/SD_{16}]])^{SD_{16}} \qquad
\text{(Proposition~\ref{prop:GHMRmapping})}. $$
\qed \end{pf*}

Neeman \cite{Neeman} defines a morphism of fiber sequences to
be \emph{good} if the induced sequence of fibers is a fiber sequence.

\begin{mylem}
The morphism of fiber sequences given by Lemma~\ref{lem:fibermap} is good.
\end{mylem}

\begin{pf*}{Proof.}
There exists a morphism $h$ such that the morphism of fiber sequences
$$
\xymatrix@C+2em{
\Sigma^{-1}Q \wedge \br{S} \ar[r]^{\delta'} \ar[d]|{1 \wedge (\tau-1)} &
DQ \wedge \br{S} \ar@{.>}[d]_{h} \ar[r]^{D\eta\wedge 1} &
S \wedge \br{S} \ar[d]|{1 \wedge(\tau-1)} \ar[r]^{\eta \wedge 1} &
Q \wedge \br{S} \ar[d]|{1 \wedge (\tau-1)} 
\\
\Sigma^{-1}Q \wedge \br{S} \ar[r]^{\delta'} &
DQ \wedge \br{S} \ar[r]^{D\eta\wedge 1} &
S \wedge \br{S} \ar[r]^{\eta \wedge 1} &
Q \wedge \br{S}
} $$
is a good morphism.  The morphism 
$$ 1 \wedge(\tau-1): DQ \wedge \br{S} \rightarrow DQ \wedge \br{S} $$
differs from $h$ by a composite $\delta' \circ \alpha \circ D\eta \wedge 1$
where $\alpha$ lies in 
$$ [\br{S}, \Sigma^{-1}Q \wedge \br{S}] \cong \pi_1 F(\br{S}, \br{Q}) \qquad
\text{(Lemma~\ref{lem:xi})}. $$
Using Theorem~\ref{thm:FEHEH}, we see that the mapping spectrum $F(\br{S},
\br{Q})$ is built from the spectra
\begin{align*}
F(\br{S}, TMF) & \simeq TMF[[\ZZ_p]] \\
F(\br{S}, \Sigma^{7}E(2)) & \simeq \Sigma^{7}E(2)[[\ZZ_p]]
\end{align*}
and both of these spectra have trivial $\pi_1$.  We conclude that 
$$ [\br{S},\Sigma^{-1}Q \wedge \br{S}] = 0 $$
and $h$ must
equal $1 \wedge(\tau-1)$.
\qed \end{pf*}

Since the morphism of fiber sequences given by Lemma~\ref{lem:fibermap} is
good, 
\emph{there exist} induced maps $f$, $g$, and $h$ 
below which make a $3\times3$
diagram of fiber sequences, as displayed below.
\begin{equation}\label{diag:main3x3}
\xymatrix{
\Sigma^{-2} Q \wedge \br{S} \ar[r] \ar@{}[dr]|{(-1)} \ar[d] &
\Sigma^{-1} DQ \wedge \br{S} \ar[d] \ar[r] &
\Sigma^{-1} \br{S} \ar[r] \ar[d] &
\Sigma^{-1} Q \wedge \br{S} \ar[d]
\\
\Sigma^{-1} Q \ar@{.>}[r]_f \ar[d] &
DQ \ar@{.>}[r]_g \ar[d] &
S \ar@{.>}[r]_h \ar[d] &
Q \ar[d]
\\
\Sigma^{-1} Q \wedge \br{S} \ar[r] \ar[d] &
DQ \wedge \br{S} \ar[d] \ar[r] &
\br{S} \ar[r] \ar[d] &
Q \wedge \br{S} \ar[d]
\\
\Sigma^{-1} Q \wedge \br{S} \ar[r] &
DQ \wedge \br{S} \ar[r] &
\br{S} \ar[r] &
Q \wedge \br{S} 
}
\end{equation}
We must show that $g = D\eta$ and $h = \eta$.  Our identifications of the
Spanier-Whitehead duals of $TMF$ and $E(2)$ in
Proposition~\ref{prop:DTMF} in particular imply that these spectra 
are reflexive.
Therefore the spectrum $\br{Q}$, being the fiber of a map between
reflexive spectra, is itself reflexive.  Given different maps $g'$ and
$h'$ in place of $g$ and $h$ making Diagram~\ref{diag:main3x3} commute, we
see that the difference $g-g'$ is in the image of the homomorphism
$$ \delta_* : [DQ,\Sigma^{-1}\br{S}] \rightarrow [DQ,S]. $$
We have
isomorphisms
\begin{alignat*}{2}
[DQ,\Sigma^{-1} \br{S} ] & \cong [DQ, D\br{S}] & \qquad &\text{(Proposition
\ref{prop:Sbar})} \\
& \cong [DQ\wedge \br{S}, S] & &\\
& \cong [\Sigma D(\br{Q}), S] & \qquad & \text{(Lemma \ref{lem:DQ^Sbar})} \\
& \cong \pi_1(DD\br{Q}) & & \\
& \cong \pi_1(\br{Q}) & \qquad & \text{($\br{Q}$ is reflexive)}
\end{alignat*}

Similarly, the difference $h-h'$ is in the image of the homomorphism
$$ (1 \wedge \delta)_* : [S,\Sigma^{-1}Q \wedge \br{S}] \rightarrow [S,Q]
$$
and we have, by Lemma~\ref{lem:xi}, an isomorphism
$$ [S,\Sigma^{-1}Q\wedge \br{S}] \cong \pi_1(\br{Q}). $$
However, $\pi_1(\br{Q}) = 0$, since $\br{Q}$ is built from $TMF$ and
$\Sigma^7 E(2)$.  Thus $g$ and $h$ are uniquely determined up to homotopy
having the property that they make Diagram~\ref{diag:main3x3} commute.
Since $D\eta$ and $\eta$ make Diagram~\ref{diag:main3x3} commute, we
must have $g = D\eta$ and $h = \eta$, as desired.  This completes the proof
of Theorem~\ref{thm:cofiber}.

\appendix
\section*{Appendix~A: A cellular model for $\br{S}$}
\addcontentsline{toc}{part}{Appendix~A: A cellular model for $\br{S}$
\hfill}
\addtocounter{section}{1}
\setcounter{subsection}{0}
\setcounter{mythm}{0}

In this appendix we shall always be working in the $K(n)$-local category at
an arbitrary prime $p$.  Let $\GG_n$ be the $n^\mathrm{th}$ 
extended Morava stabilizer
group, and let $\GG^1_n$ be the kernel of the reduced norm
\cite{DevinatzHopkins}
$$ 1 \rightarrow  \GG^1_n \rightarrow \GG_n \rightarrow \ZZ_p \rightarrow 1. $$
Let $\br{S}$ be the homotopy fixed point spectrum $E_n^{h\GG^1}$.
Let $\tau \in \ZZ_p$ be a topological generator.
In \cite{DevinatzHopkins}, the following theorem is proven as an
application of the authors' continuous homotopy fixed point construction.

\begin{mythm}[Hopkins-Miller]
There is a cofiber sequence
$$ \Sigma^{-1}\br{S} \xrightarrow{\delta} S \xrightarrow{r} \br{S}
\xrightarrow{[\tau]-1} \br{S}. $$
The element $\zeta$ exists in $\pi_{-1}(S)$, and is given by the
composite
$$ \zeta: S^{-1} \xrightarrow{r} \Sigma^{-1}\br{S}
\xrightarrow{\delta} S. $$
\end{mythm}

We shall prove the following theorem, which gives a $K(n)$-local cellular
decomposition of $\br{S}$.  We used this decomposition in the proof of 
Lemma~\ref{lem:xi} to
compute $V(1)_*(Q(2)\wedge \br{S})$.

\begin{mythm}\label{thm:Sbarcell}
There exist complexes $C^i_\zeta$ with cellular decompositions
$$ C^i_\zeta = S^0 \cup_\zeta \underbrace{e^0 \cup_\zeta e^0 \cup_\zeta \cdots
\cup_\zeta e^0}_{i}. $$
There is an equivalence $\varinjlim C^i_\zeta \simeq \br{S}$.
\end{mythm}

\begin{myrem}
The existence of the complexes $C^i_\zeta$ is equivalent to the the Toda
bracket
$$ \bra{\underbrace{\zeta, \ldots, \zeta}_{i}}$$
being defined and containing zero, for every $i$.
\end{myrem}

The remainder of this appendix is dedicated to proving this theorem.
Our models for the intermediate complexes $C^i_\zeta$ will be the 
homotopy fibers of the map $([\tau]-1)^i$, giving fiber sequences
$$ \Sigma^{-1} \br{S} \xrightarrow{\delta_i} 
C^i_\zeta \xrightarrow{r_i} \br{S} \xrightarrow{([\tau]-1)^i} \br{S}. $$
The complex $C^0_\zeta$ is just the sphere $S$.  

We must explain why the homotopy fibers $C^i_\zeta$ have the cellular
models as claimed. 
We shall inductively prove that there are cofiber sequences
\begin{equation}\label{eq:Czetacellular}
S^{-1} \xrightarrow{\zeta_{i-1}} C^{i-1}_\zeta \xrightarrow{\iota_i} 
C^i_\zeta \xrightarrow{\nu_i} S
\end{equation}
where the map $\zeta_{i-1}$ will make the following diagram commute.
$$
\xymatrix{
S^{-1} \ar[r]^{\zeta_{i-1}} \ar[dr]_\zeta &
C^{i-1}_\zeta \ar[d]^{\nu_{i-1}} 
\\
& S}
$$
Thus $\zeta_{i-1}$ is a lift of $\zeta$ to the whole complex
$C^{i-1}_\zeta$.

The existence of the cofiber sequence~\ref{eq:Czetacellular} is a direct
consequence of Verdier's axiom.
\begin{equation}\label{diag:3x3zeta} 
\xymatrix@C+1em@R+1em{
\Sigma^{-1}\br{S} \ar@{=}[r] \ar[d]_{([\tau]-1)^{i-1}} &
\Sigma^{-1}\br{S} \ar[r] \ar[d]|{([\tau]-1)^{i}} &
\ast \ar[r] \ar[d] &
\br{S} \ar[d]^{([\tau]-1)^{i-1}}
\\
\Sigma^{-1}\br{S} \ar[r]^{[\tau]-1} \ar[d]_{\delta_{i-1}} &
\Sigma^{-1}\br{S} \ar[d]_{\delta_i} \ar[r]_\delta &
S \ar[r]_r \ar@{=}[d] &
\br{S} \ar[d]^{\delta_{i-1}}
\\
C^{i-1}_\zeta \ar@{.>}[r]_{\iota_i} \ar[d]_{r_{i-1}} &
C^i_\zeta \ar[d]_{r_i} \ar@{.>}[r]_{\nu_i} &
S \ar[r]_{\zeta_{i-1}} \ar[d] &
\Sigma C^{i-1}_\zeta \ar[d]^{r_{i-1}}
\\
\br{S} \ar@{=}[r] &
\br{S} \ar[r] &
\ast \ar[r] &
\Sigma \br{S} 
}
\end{equation}

The commutativity of the following diagram implies that
the composite of $\zeta_{i}$ with the projection onto the top cell
is $\zeta$.
$$
\xymatrix{
S^{-1} \ar[rr]^{\zeta_{i}} \ar[rd]^r \ar@/_2pc/[ddrr]_\zeta & &
C^{i}_\zeta \ar[dd]^{\nu_i} 
\\
& \Sigma^{-1} \br{S} \ar[dr]^\delta \ar[ur]^{\delta_i} 
\\
& & S
} $$

In Diagram~$\ref{diag:3x3zeta}$, we saw that there are commutative diagrams
$$ 
\xymatrix{
C^{i-1}_\zeta \ar[rr]^{\iota_i} \ar[dr]_{r_{i-1}} & & C^i_\zeta \ar[dl]^{r_i}
\\
& \br{S}
} $$
Let $C^\infty_\zeta$ denote the homotopy colimit $\varinjlim C^i_\zeta$.
The compatibility of the $r_i$'s implies that there is a map
$$ r_\infty = \varinjlim r_i : C^\infty_\zeta \rightarrow \br{S}. $$

We are left with showing that
the map $r_\infty$ is an equivalence.  
The left hand columns of Diagram~\ref{diag:3x3zeta} give us a tower of
cofiber sequences
$$
\xymatrix{
\Sigma^{-1}\br{S} \ar[r]^{\delta_{i-1}} \ar[d]^{[\tau]-1} &
C^{i-1}_\zeta \ar[r]^{r_{i-1}} \ar[d]^{\iota_{i}} &
\br{S} \ar@{=}[d] 
\\
\Sigma^{-1}\br{S} \ar[r]_{\delta_i} &
C^i_\zeta \ar[r]_{r_i} &
\br{S}
}
$$
which, upon taking homotopy colimits gives a cofiber sequence
$$ \Sigma^{-1}([\tau]-1)^{-1}\br{S} \xrightarrow{\delta_{\infty}} 
C^\infty_\zeta \xrightarrow{r_\infty} \br{S}. $$
Therefore, it suffices to prove the following lemma.  The author is 
grateful to
Daniel Davis for helping to streamline the proof of this lemma.

\begin{mylem}\label{lem:telnull}
The telescope $([\tau]-1)^{-1}\br{S}$ is contractible.
\end{mylem}

\begin{pf*}{Proof.}
Let $w$ denote the self-map $[\tau]-1$.
It suffices to show that $(K_n)_*(w^{-1}\br{S})$ is zero, where $K_n =
E_n/(p,u_1, \ldots, u_{n-1})$ is
$2$-periodic Morava $K$-theory.  
The $E_n$-homology of
$\br{S}$ is given by
$$ (E_n)_*(\br{S}) \cong \map^c(\ZZ_p, (E_n)_*). $$
In particular, $(E_n)_*(\br{S})$ is pro-free over $(E_n)_*$ by Theorem~2.5
of \cite{HoveyOp}. By Corollary~5.2 of \cite{HoveySS} we have an
isomorphism
$$ (K_n)_*(\br{S}) \cong \map^c(\ZZ_p, (K_n)_*). $$
Using the fact that $(K_n)_*$ is a discrete module, we have
\begin{align*}
(K_n)_*(w^{-1}\br{S}) & = \varinjlim_{w_*} (K_n)_*(\br{S}) \\
& = \varinjlim_{w_*} \map^c(\ZZ_p, (K_n)_*) \\
& = \varinjlim_{w_*} \varinjlim_{k} \map(\ZZ/p^k, (K_n)_*) \\
& = \varinjlim_{k} \varinjlim_{w_*} \map(\ZZ/p^k, (K_n)_*).
\end{align*}
We are reduced to showing that
$$ \varinjlim_{w_*} \map(\ZZ/p^k, (K_n)_*) = 0. $$
Using the isomorphism
$$ \map(\ZZ/p^k, (K_n)_*) \cong \Hom_{\FF_p}(\FF_p[\ZZ/p^k], (K_n)_*), $$
the map $w_*$ acts by right multiplication by 
$[\br{\tau}]-1$ on the $\FF_p[\ZZ/p^k]$ factor.  
Here $\br{\tau}$ is a generator of the
group $\ZZ/p^k$.  We have, since we are now in characteristic $p$, 
$$ ([\br{\tau}]-1)^{p^k} = [\br{\tau}]^{p^k} - 1 = 1-1 = 0. $$
Thus on $\map(\ZZ/p^k, (K_n)_*)$, we have $w_*^{p^k} = 0$, and the colimit
over $w_*$ is zero.
\qed \end{pf*}

\end{document}